\documentclass[10pt]{article}

\usepackage{a4wide}
\usepackage{amssymb}
\usepackage{amsfonts}
\usepackage{amsmath}
\input xy
\xyoption{arrow} \xyoption{matrix}

\date{}

\newtheorem{proposition}{Proposition}[section]
\newtheorem{theorem}[proposition]{Theorem}
\newtheorem{lemma}[proposition]{Lemma}

\newtheorem{corollary}[proposition]{Corollary}

\def\der{\partial }

\def\nFM0{{\nu }_{F,M_0}}
\def\nFN0{{\nu }_{F,N_0}}
\def\nGN0{{\nu }_{G,N_0}}

\def\N0{ {\bf N}_0 }

\def\v{\varphi}
\def\ra{\rightarrow}

\def\Xpm{X^{\pm }}

\def\s{\sigma}
\def\Z{\mathbb{Z}}

\def\l1{{\lambda}_1}

\def\a{\alpha}
\def\a0{ {\alpha }_0}
\def\a1{ {\alpha }_1}

\def\l{\lambda}
\def\o{\omega}

\def\nFGM0{{\nu }_{F,G,M_0}}

\def\nFN0{{\nu}_{F,N_0}}

\def\sm{{\sigma}^m}

\def\sm1{{\sigma}^{-1}}

\def\smtp1{{\sigma}^{-t+1}}

\def\o{\omega }
\def\S1{S^{-1}}

\def\Xpm1{X^{\pm 1}_1}

\def\sPM1{{\sigma }^{\pm 1}}
\def\sMP1{{\sigma }^{\mp 1 }}

\def\b{\beta}

\def\di{{\rm d.ind}}

\def\L{\Lambda}

\def\Ytm1{Y^{t-1}}
\def\Yim1{Y^{i-1}}

\def\CL{{\cal L}}

\def\CN{{\cal N}}

\def\ass{{\rm ass}}

\def\Aut{{\rm Aut}}

\def\bA{\overline{A}}

\def\ker{ {\rm ker } }

\def\gr{ {\rm gr} }

\def\SL2Z{ {\rm SL}_2({\bf Z}) }

\def\CR{ {\cal R}}

\def\CL{{\cal L}}

\def\Gp1{ G^{1 , 1 } }
\def\P11{ P^{-1 , 1 } }
\def\Pp1{ P^{1 , 1 } }

\def\nCLsr{{}^\nu\kern-2pt {\cal L}^{\sigma , \rho  }}
\def\nP{{}^\nu \kern-2pt P}
\def\nL{{}^\nu\kern-2pt L}
\def\nLL{{}^\nu\kern-2pt \Lambda}
\def\nPsr{{}^\nu\kern-2pt P^{\sigma , \rho  }}
\def\nLsr{{}^\nu\kern-2pt L^{\sigma , \rho  }}
\def\nuCL{{}^\nu\kern-2pt  {\cal L}}
\def\nCLsr{{}^\nu\kern-2pt {\cal L}^{\sigma , \rho  }}
\def\nCL1m{{}^\nu\kern-2pt {\cal L}^{-1 , 1  }}
\def\x1nu{x^\frac{1}{\nu}}
\def\xm1nu{x^{-\frac{1}{\nu}}}

\def\rad{{\rm rad}}

\def\ob{\overline{b}}

\def\CR{ {\cal R}}

\def\CN{{\cal N}}
\def\ra{\rightarrow }

\def\CB{{\cal B}}

\def\CI{{\cal I}}

\def\CT{{\cal T}}

\def\CC{ {\cal C}}

\def\nAM0{{\nu }_{{\cal A},M_0}}
\def\nAN0{{\nu }_{{\cal A},N_0}}

\def\CR{ {\cal R }}

\def\bR{\overline{R}}

\def\ga{\mathfrak{a}}
\def\gb{\mathfrak{b}}
\def\gc{\mathfrak{c}}

\def\gn{\mathfrak{n}}

\def\gr{\mathfrak{r}}

\def\SL{{\rm SL}}

\def\di!{\frac{\der^i}{i!}}
\def\dik!{\frac{\der^k_i}{k!}}

\def\Max{{\rm Max}}

\def\N{\mathbb{N}}

\def\0{\overline{0}}
\def\1{\overline{1}}

\def\Ln1{\L_{n,\overline{1}}}

\def\a1{a_{\overline{1}}}

\def\bs{\overline{s}}

\def\S{\Sigma}

\def\vn1{\overrightarrow{n-1}}

\def\Q{\mathbb{Q}}

\def\Min{{\rm Min}}

\def\Inn{{\rm Inn}}

\def\mJ{\mathbb{J}}
\def\mI{\mathbb{I}}

\def\K1{{\rm K}_1}

\def\hmI1{\widehat{\mI_1}}
\def\tmI1{\widetilde{\mI_1}}
\def\tmJ1{\widetilde{\mJ_1}}
\def\hB1{\widehat{B_1}}
\def\hCB1{\widehat{\CB_1}}

\def\bS{\overline{S}}

\def\Den{{\rm Den}}
\def\Denl{{\rm Den}_l}
\def\Ore{{\rm Ore}}

\def\Den{{\rm Den}}

\def\Loc{{\rm Loc}}

\def\Ass{{\rm Ass}}

\def\maxDen{{\rm max.Den}}
\def\maxAss{{\rm max.Ass}}
\def\maxLoc{{\rm max.Loc}}
\def\llrad{{\rm l.lrad}}

\def\assmaxDen{{\rm ass.max.Den}}

\def\br{\overline{r}}

\def\bs{\overline{s}}

\def\ga{\mathfrak{a}}

\def\tor{{\rm tor}}

\def\bE{\overline{E}}

\def\gll{\mathfrak{l}}

\def\IDen{{\rm IDen}}
\def\ILoc{{\rm ILoc}}
\def\maxker{{\rm max.ker}}
\def\CII{{\cal II}}
\def\CIT{{\cal IT}}

\begin{document}

\author{V. V. \  Bavula 
}

\title{Left localizations of left Artinian rings}

\maketitle

\begin{abstract}

For an arbitrary left Artinian ring $R$, explicit descriptions  are given of all the left  denominator sets $S$ of $R$ and   left localizations $S^{-1}R$ of $R$. It is proved that,  up to $R$-isomorphism, there are only finitely many  left localizations and  each of them  is an idempotent localization, i.e. $S^{-1}R\simeq S_e^{-1}R$ and $\ass (S) = \ass (S_e)$ where $S_e=\{1,e\}$ is a left denominator set of $R$  and $e$ is an idempotent. Moreover, the idempotent $e$  is unique up to a conjugation. It is proved that  the number of maximal left denominator sets of $R$ is finite and does not exceed the number of isomorphism classes of simple left $R$-modules. The set of maximal left denominator sets of $R$  and the  left localization radical of $R$ are described.

$\noindent $



$${\bf Contents}$$
\begin{enumerate}
\item Introduction.
\item Preliminaries.
\item Idempotent left denominator sets.
 \item Left localizations of left Artinian rings.
\item Structure of left Artinian rings with zero left localization radical.
 \item Characterization of the left localization radical of a left Artinian ring.
 \item Description of left denominator sets of a left Artinian ring.
 \item Localizations of Artinian rings.
\item  Rings with  left Artinian left quotient ring.
\end{enumerate}
\end{abstract}


\section{Introduction}

In this paper, module means a left module, and the following notation is fixed:
\begin{itemize}
\item  $R$ is a ring with 1,  $R^*$ is its group of units and $\Inn (R):=\{ \o_u\, | \,  u\in R^*\}$  is the group of inner automorphisms of $R$ where $\o_u(r):= uru^{-1}$ for $r\in R$, $\rad (R)$ is the Jacobson radical of $R$;
\item $\Ore_l(R):=\{ S\, | \, S$ is a left Ore set in $R\}$; \item
$\Den_l(R):=\{ S\, | \, S$ is a left denominator set in $R\}$;
\item $\Loc_l(R):= \{ [S^{-1}R]\, | \, S\in \Den_l(R)\}$ is the set of $R$-{\em isomorphism  classes} $[S^{-1}R]$ of left localizations of the ring $R$;  $[S^{-1}R]=[S'^{-1}R]$ iff
 the map $S^{-1}R \ra S'^{-1}R$, $s^{-1}r\mapsto s^{-1}r$ is a well-defined isomorphism. We identify  $S^{-1}R$  with $[S^{-1}R]$. So $S^{-1}R = S'^{-1}R$
  iff the rings $S^{-1}R$ and $ S'^{-1}R$ are $R$-isomorphic;
   \item
$\Ass_l(R):= \{ \ass (S)\, | \, S\in \Den_l(R)\}$ where $\ass
(S):= \{ r\in R \, | \, sr=0$ for some $s=s(r)\in S\}$;
\item $\Den_l(R, \ga ) :=\{ S\in \Den_l(R) \, | \, \ass (S) = \ga \}$.
\end{itemize}

In brief, for an arbitrary left Artinian  $R$, this paper presents a complete picture of how left localizations and left denominator sets of $R$ look like, and the situation is so `simple' and natural that one cannot image a better/simpler one, see below.

$\noindent $

{\bf Every left localization of left Artinian ring  is an idempotent left localization}.   We say that a left localization $A=S^{-1}R$  of a ring $R$ is an {\em idempotent left localization} of $R$ if there is an idempotent $e\in R$ such that $S_e= \{ 1,e\}\in \Den_l(R)$ and the rings   $S_e^{-1}R$ and $ A$ are $R$-{\em isomorphic} (equivalently, the map $A=S^{-1}R\ra S_e^{-1}R$, $s^{-1}r\mapsto s^{-1}r$, is an isomorphism). 

The  following theorem shows that every left localization of a left Artinian ring is an idempotent left localization.

\begin{itemize}
\item {\bf (Theorem \ref{A9Feb13})}
{\em Let $R$ be a left Artinian ring and $S\in \Den_l(R, \ga )$. Then
 there exists a nonzero idempotent $e\in R$ such that $S_e:=\{ 1,e\}\in \Den_l(R, \ga )$ and  the rings  $S^{-1}R$ and $ S_e^{-1}R$ are $R$-isomorphic. The idempotent $e$ is unique up to conjugation.}
\end{itemize}

{\bf There are only finitely many left localizations for a left Artinian ring}. Let $R$ be a left Artinian ring, $\rad (R)$ be its radical, $\bR:= R/ \rad (R)=\prod_{i=1}^s\bR_i$ -- a direct  product of simple Artinian rings $\bR_i$, $\overline{1}_i$ be the identity element of the ring $\bR_i$. So, $1= \sum_{i=1}^s \overline{1}_i$ is the sum of orthogonal central idempotents of $\bR$, $1=\sum_{i=1}^s 1_i$ is a sum of orthogonal  idempotents of $R$ such that $1_i$ is a lifting of $\overline{1}_i$ (see (\ref{1=S1is})). For each non-empty set $I$ of $\{ 1, \ldots , s\}$, let $e_I:=\sum_{i\in I}1_i$,
$$ \CI_l':= \CI_l'(R):=\{ e_I\, | \, e_IR(1-e_I)=0\}\;\; {\rm and}\;\; | \CI_l'|<\infty .$$

 The following theorem shows that for a left Artinian ring $R$  there are only finitely many left localizations. Moreover, it gives  explicit descriptions of the sets $\Loc_l(R)$ and $\Ass_l(R)$.
\begin{itemize}
\item {\bf (Theorem \ref{A11Feb13})}
{\em Let $R$ be a left Artinian ring. Then
 the map $\CI_l'(R)\ra \Loc_l(R)$, $e\mapsto S_e^{-1}R=R/(1-e)R$, is a bijection. The map $\CI_l'(R)\ra \Ass_l(R)$, $e\mapsto (1-e)R$, is a bijection.
}
\end{itemize}

{\bf Classification of  denominator sets of a left Artinian ring}. A subset $S$ of a ring $R$ is called a {\em multiplicative set} if $1\in S$, $SS\subseteq S$ and $0\not\in S$.  The next theorem gives  a criterion for a multiplicative set of a left Artinian ring to be a left denominator set and provides  an explicit  description/classification of all the left denominator sets of $R$.

\begin{itemize}
\item {\bf (Theorem \ref{F12Feb13})}
{\em Let $R$ be a left Artinian ring and $S$ be a multiplicative set of $R$. The following statements are equivalent.}
\begin{enumerate}
\item $S\in \Den_l(R)$.
\item {\em There is a nonzero idempotent $e\in R$ such that $eR(1-e)=0$, $S\subseteq \begin{pmatrix}
 R_{11}^* & 0\\ R_{21} & R_{22}
 \end{pmatrix}$ and there is an element $s\in S$ such that $s= \begin{pmatrix}
 u& 0\\ v & 0
 \end{pmatrix}$ where $R=\begin{pmatrix}
 R_{11}& 0\\ R_{21} & R_{22}
 \end{pmatrix}$ is the matrix ring associated with the idempotent $e$ (see (\ref{MRawe})).
\item There is a unit $\l \in R^*$, and idempotent $e\in \CI'_l(R)$ and an element $s\in S$ such that $\l S\l^{-1} \subseteq \begin{pmatrix}
 R_{11}^*& 0\\ R_{21} & R_{22}
 \end{pmatrix}$ and $\l s\l^{-1} = \begin{pmatrix}
 u& 0\\ v & 0
 \end{pmatrix}$ where $R=\begin{pmatrix}
 R_{11}& 0\\ R_{21} & R_{22}
 \end{pmatrix}$ is the matrix ring associated with the idempotent $e$.
 \item There is an element $s\in S$ such that $S_s\in \Den_l(R)$ where $S_s\;= \{ s^i\,  | \, i\in \N\}$  and the images of all the elements of $S$ in the ring $S_s^{-1}R$ are units.}
\end{enumerate}
{\em If one of the equivalent conditions holds then $\ass (S) = \ker_R(s\cdot  )$ in all three cases regardless of the choice of $s$.}
\end{itemize}

So, in order to obtain all the left denominator sets of a left Artinian ring  $R$ we have to choose an element $e\in \CI_l'(R)$, a multiplicative set $S$ of $\begin{pmatrix}
 R_{11}^*& 0\\ R_{21} & R_{22}
 \end{pmatrix}$ that contains an element of the form $\begin{pmatrix}
 u& 0\\ v & 0
 \end{pmatrix}$. Then $S$  is a left denominator set of $R$ and an arbitrary left denominator set of $R$ is of the type $\l S\l^{-1}$ for some $\l\in R^*$.

{\bf The maximal left denominator sets of a left Artinian ring $R$}. In \cite{larglquot}, the concept of maximal left denominator set of ring was introduced and it was shown that the set $\maxDen_l(R)$ of   maximal left denominator sets of a ring $R$ is a non-empty set. For a left Artinian ring $R$,   the finite set $(\CI_l', \geq )$ is a partially ordered set where $e_I\geq e_J$ iff $I\supseteq J$. Let $\min \CI_l'$ be the set of minimal elements of $\CI'$. The next theorem provides a description of the maximal left denominator sets of a left Artinian ring.

\begin{itemize}
\item {\bf (Theorem \ref{B11Feb13})}
{\em Let $R$ be a left Artinian ring. Then
 \begin{enumerate}
\item $\maxDen_l(R)=\{ T_e\, | \, e\in \min \CI_l'(R)\} $ where $T_e=\{ u\in R\, | \, u+(1-e)R\in (R/(1-e)R)^*\}$.
\item $|\maxDen_l(R)|\leq s$  where $s$ is  the number of isomorphism classes of left simple $R$-modules.
    \item  $|\maxDen_l(R)|= s$ iff $R$ is a semisimple ring.
\end{enumerate}
  }
\end{itemize}
{\bf The maximal left denominator sets of a ring with left Artinian left quotient ring}. The next theorem shows that a ring  with left Artinian left quotient ring has only finitely many  maximal left denominator sets.
\begin{itemize}
\item {\bf (Theorem \ref{X17Feb14})}
{\em Let $R$ be a ring such that  $Q_l(R)$ is a left Artinian ring and $s$ be the number of iso-classes of simple  left $Q_l(R)$-modules. Then}
\begin{enumerate}
\item {\em the map $\maxDen_l(R)\ra \maxDen_l(Q_l(R))$, $S\mapsto SQ_l(R)^*$, is a bijection with the inverse $T\mapsto T\cap R$. In particular, $|\maxDen_l(R)|=|\maxDen_l(Q_l(R))|\leq s<\infty$.
    \item $|\maxDen_l(R)|=s$ iff $Q_l(R)$ is a semisimple ring iff $R$ is a semiprime left Goldie ring.}
\end{enumerate}
\end{itemize}
 Recall that the largest left quotient ring $Q_l(R)$  of $R$ is a left Artinian ring   iff the (classical)  left quotient ring $Q_{l, cl} (R):= \CC_R^{-1}R$ is a left Artinian ring, and in this case $S_l(R) = \CC_R$, \cite{Bav-genGoldie}.

{\bf Criterion for the powers of an element to be a left denominator set}. For a left Artinian ring $R$, the following theorem is an explicit criterion for the powers of a non-nilpotent element of $R$ to be a left denominator set.
\begin{itemize}
\item {\bf (Theorem \ref{C11Feb13})}
{\em Let $R$ be a left Artinian ring, $\CI_l'(R)$ be as above, $s\in R$ be a non-nilpotent element of $R$, $e=e(s)$ be the idempotent associated with the element $s$ (see (\ref{1=ese})), $S_e=\{ 1,e\}$ and $S_s=\{ s^i\, | \, i\in \N\}$. The following statements are equivalent.
\begin{enumerate}
\item $S_s\in \Den_l(R)$.
\item $S_e\in \Den_l(R)$ and $(1-e)s(1-e)$ is a nilpotent element.
\item $eR(1-e)=0$ and $(1-e)s(1-e)$ is a nilpotent element.
\end{enumerate}
If one of the equivalent conditions holds then $\ass (S_s) = (1-e)R$ and $S_s^{-1}R \simeq S_e^{-1}R \simeq R/ (1-e)R$, the core $S_{s,c}$ of the left denominator set $S_s$ is equal to $\{ s^i\, | \, i\geq 1, (1-e) s^i(1-e)=0\}$.
  }
\end{itemize}
{\bf Duality between left and right localizations of an Artinian ring}. The sets of left and right localizations of a ring $R$, $(\Loc_l(R), \ra )$ and $(\Loc_r(R), \ra )$, are partially ordered sets where $[A]\ra [B]$ if there is a ring  $R$-homomorphism  $A\ra B$. In general, the left and right localizations of a ring $R$ are almost unrelated but for each Artinian ring $R$ there is a duality between the partially ordered sets
 $(\Loc_l(R), \ra )$ and $(\Loc_r(R), \ra )$.
 \begin{itemize}
 \item {\bf (Theorem \ref{25Apr14})}
{\em Let $R$ be an Artinian ring. Then the map
$$\Loc_l(R)\backslash [R]\ra \Loc_r(R)\backslash [R], \;\; [ R/ (1-e_I)R]\ra [ R/ R(1-e_{CI})], $$ is an anti-isomorphism of posets (i.e. an order reversing bijection). In particular, $|\Loc_l(R)|= |\Loc_r(R)|$. }
\end{itemize}


\section{Preliminaries}\label{PRLM}

In this section, we collect necessary results that are used in the proofs of this paper. More results on localizations of rings (and some of the missed standard definitions) the reader can find in \cite{Jategaonkar-LocNRings}, \cite{Stenstrom-RingQuot} and \cite{MR}.  In this paper the following notation will remained fixed.

$\noindent $

{\bf Notation}:

\begin{itemize}

\item $S_\ga=S_\ga (R)=S_{l,\ga }(R)$
 is the {\em largest element} of the poset $(\Den_l(R, \ga ),
\subseteq )$ and $Q_\ga (R):=Q_{l,\ga }(R):=S_\ga^{-1} R$ is  the
{\em largest left quotient ring associated to} $\ga$, $S_\ga $
exists (Theorem 
 2.1, \cite{larglquot});
\item in particular, $S_0=S_0(R)=S_{l,0}(R)$ is the largest
element of the poset $(\Den_l(R, 0), \subseteq )$ and
$Q_l(R):=S_0^{-1}R$ is the {\em largest left quotient ring} of $R$, \cite{larglquot};
\item $\Loc_l(R, \ga ):= \{ S^{-1}R\in \Loc_l(R)\, | \, S\in \Den_l(R, \ga
)\}$.
\end{itemize}

In \cite{larglquot}, we introduce the following new concepts and
prove their existence for an {\em arbitrary} ring: {\em the
largest left quotient ring of a ring, the largest  regular left
 Ore  set of a ring, the left localization radical of a ring, a maximal left denominator set,  a maximal left quotient ring of a ring,
 a (left) localization
maximal ring}.
Using an analogy with rings, the counter parts of the three concepts: a maximal left denominator set, the left localization radical and a maximal left quotient ring,
 for rings would be a left maximal ideal, the Jacobson radical and  a simple factor ring,  respectively.

{\bf The largest regular left Ore set and the largest left
quotient ring of a ring}. Let $R$ be a ring. A {\em
multiplicatively closed subset} $S$ of $R$ or a {\em
 multiplicative subset} of $R$ (i.e. a multiplicative sub-semigroup of $(R,
\cdot )$ such that $1\in S$ and $0\not\in S$) is said to be a {\em
left Ore set} if it satisfies the {\em left Ore condition}: for
each $r\in R$ and
 $s\in S$, $ Sr\bigcap Rs\neq \emptyset $.
Let $\Ore_l(R)$ be the set of all left Ore sets of $R$.
  For  $S\in \Ore_l(R)$, $\ass (S) :=\{ r\in
R\, | \, sr=0 \;\; {\rm for\;  some}\;\; s\in S\}$  is an ideal of
the ring $R$.

$\noindent $

A left Ore set $S$ is called a {\em left denominator set} of the
ring $R$ if $rs=0$ for some elements $ r\in R$ and $s\in S$ implies
$tr=0$ for some element $t\in S$, i.e. $r\in \ass (S)$. Let
$\Den_l(R)$ be the set of all left denominator sets of $R$. For
$S\in \Den_l(R)$, let $S^{-1}R=\{ s^{-1}r\, | \, s\in S, r\in R\}$
be the {\em left localization} of the ring $R$ at $S$ (the {\em
left quotient ring} of $R$ at $S$). Let us stress that in Ore's method of localization one can localize {\em precisely} at left denominator sets.

In general, the set $\CC$ of regular elements of a ring $R$ is
neither left nor right Ore set of the ring $R$ and as a
 result neither left nor right classical  quotient ring ($Q_{l,cl}(R):=\CC^{-1}R$ and
 $Q_{r,cl}(R):=R\CC^{-1}$) exists.
 Remarkably, there  exists the largest
 regular left Ore set $S_0= S_{l,0} = S_{l,0}(R)$, \cite{larglquot}. This means that the set $S_{l,0}(R)$ is an Ore set of
 the ring $R$ that consists
 of regular elements (i.e., $S_{l,0}(R)\subseteq \CC$) and contains all the left Ore sets in $R$ that consist of
 regular elements. Also, there exists the largest regular (left and right) Ore set $S_{l,r,0}(R)$ of the ring $R$.
 In general, all the sets $\CC$, $S_{l,0}(R)$, $S_{r,0}(R)$ and $S_{l,r,0}(R)$ are distinct, for example,
 when $R= \mI_1= K\langle x, \der , \int\rangle$  is the ring of polynomial integro-differential operators  over a field $K$ of characteristic zero,  \cite{Bav-intdifline}. In  \cite{Bav-intdifline},  these four sets are found for $R=\mI_1$.

$\noindent $

{\it Definition}, \cite{Bav-intdifline}, \cite{larglquot}.    The ring
$$Q_l(R):= S_{l,0}(R)^{-1}R$$ (respectively, $Q_r(R):=RS_{r,0}(R)^{-1}$ and
$Q(R):= S_{l,r,0}(R)^{-1}R\simeq RS_{l,r,0}(R)^{-1}$) is  called
the {\em largest left} (respectively, {\em right and two-sided})
{\em quotient ring} of the ring $R$.

$\noindent $

 In general, the rings $Q_l(R)$, $Q_r(R)$ and $Q(R)$
are not isomorphic, for example, when $R= \mI_1$, \cite{Bav-intdifline}.  The next
theorem gives various properties of the ring $Q_l(R)$. In
particular, it describes its group of units.


\begin{theorem}\label{4Jul10}
\cite{larglquot}
\begin{enumerate}
\item $ S_0 (Q_l(R))= Q_l(R)^*$ {\em and} $S_0(Q_l(R))\cap R=
S_0(R)$.
 \item $Q_l(R)^*= \langle S_0(R), S_0(R)^{-1}\rangle$, {\em i.e. the
 group of units of the ring $Q_l(R)$ is generated by the sets
 $S_0(R)$ and} $S_0(R)^{-1}:= \{ s^{-1} \, | \, s\in S_0(R)\}$.
 \item $Q_l(R)^* = \{ s^{-1}t\, | \, s,t\in S_0(R)\}$.
 \item $Q_l(Q_l(R))=Q_l(R)$.
\end{enumerate}
\end{theorem}

{\bf The maximal left denominator sets and the maximal left localizations  of a ring}. The set $(\Den_l(R), \subseteq )$ is a poset (partially ordered
set). In \cite{larglquot}, it is proved  that the set
$\maxDen_l(R)$ of its maximal elements is a {\em non-empty} set.

$\noindent $

{\it Definition}, \cite{larglquot}. An element $S$ of the set
$\maxDen_l(R)$ is called a {\em maximal left denominator set} of
the ring $R$ and the ring $S^{-1}R$ is called a {\em maximal left
quotient ring} of the ring $R$ or a {\em maximal left localization
ring} of the ring $R$. The intersection
\begin{equation}\label{llradR}
\gll_R:=\llrad (R) := \bigcap_{S\in \maxDen_l(R)} \ass (S)
\end{equation}
is called the {\em left localization radical } of the ring $R$,
\cite{larglquot}.

$\noindent $

 For a ring $R$, there is the canonical exact
sequence 
\begin{equation}\label{llRseq}
0\ra \gll_R \ra R\stackrel{\s }{\ra} \prod_{S\in \maxDen_l(R)}S^{-1}R, \;\; \s := \prod_{S\in \maxDen_l(R)}\, \s_S,
\end{equation}
where $\s_S:R\ra S^{-1}R$, $r\mapsto \frac{r}{1}$. For a ring $R$ with a semisimple left quotient ring, the left localization radical $\gll_R$ coincides with the prime radical $\gn_R$ of $R$, \cite{larglquot}. In general, $ \gll_R\neq \gn_R$ and $\gll_R\neq \rad (R)$, Theorem \ref{D12Feb13}.(4).

$\noindent $

{\bf The maximal elements of $\Ass_l(R)$}.  Let $\maxAss_l(R)$ be
the set of maximal elements of the poset $(\Ass_l(R), \subseteq )$
and

\begin{equation}\label{mADen}
\assmaxDen_l(R) := \{ \ass (S) \, | \, S\in \maxDen_l(R) \}.
\end{equation}
These two sets are equal (Proposition \ref{b27Nov12}), a proof is
based on Lemma \ref{1a27Nov12} and Corollary \ref{d4Jan13}. For an non-empty set $X$ or $R$, let ${\rm  r.ass} (X):=\{ r\in R\, | \, rx=0$ for some $x=x(r)\in X \}$.

\begin{lemma}\label{1a27Nov12}
 \cite{larglquot}
Let $S\in \Den_l(R, \ga )$ and $T\in \Den_l(R, \gb )$ be such that $ \ga \subseteq \gb$. Let $ST$ be the multiplicative semigroup generated by $S$ and $T$ in $(R,\cdot )$.  Then
\begin{enumerate}
\item ${\rm  r.ass} (ST)\subseteq \gb$.
\item $ST \in \Den_l(R, \gc )$ and $\gb \subseteq \gc$.
\end{enumerate}
\end{lemma}

\begin{corollary}\label{d4Jan13}
Let $R$ be a ring, $S\in \maxDen_l(R)$ and $T\in \Den_l(R)$. Then $T\subseteq S$ iff $\ass (T)\subseteq \ass (S)$.

\end{corollary}

{\it Proof}. $(\Rightarrow )$ If $T\subseteq S$ then $\ass (T)\subseteq \ass (S)$.

$(\Leftarrow )$ If $\ass (T)\subseteq \ass (S)$. then, by Lemma \ref{1a27Nov12}, $ST\in \Den_l(R)$ and $S\subseteq ST$, hence $S= ST$, by the maximality of $S$. Then $T\subseteq S$.  $\Box $


\begin{proposition}\label{b27Nov12}
\cite{larglquot} $\; \maxAss_l(R)= \assmaxDen_l(R)\neq \emptyset$. In particular, the ideals of this set are incomparable (i.e. neither $\ga\nsubseteq \gb$ nor $\ga\nsupseteq \gb$).
\end{proposition}

{\bf Properties of the maximal left quotient rings of a ring}.
The next theorem describes various properties of the maximal left
quotient rings of a ring, in particular, their groups of units and
their largest left quotient rings. 

\begin{theorem}\label{15Nov10}
\cite{larglquot} Let $S\in \maxDen_l(R)$, $A= S^{-1}R$, $A^*$ be
the group of units of the ring $A$; $\ga := \ass (S)$, $\pi_\ga
:R\ra R/ \ga $, $ a\mapsto a+\ga$, and $\s_\ga : R\ra A$, $
r\mapsto \frac{r}{1}$. Then
\begin{enumerate}
\item $S=S_\ga (R)$, $S= \pi_\ga^{-1} (S_0(R/\ga ))$, $ \pi_\ga
(S) = S_0(R/ \ga )$ and $A= S_0( R/\ga )^{-1} R/ \ga = Q_l(R/ \ga
)$. \item  $S_0(A) = A^*$ and $S_0(A) \cap (R/ \ga )= S_0( R/ \ga
)$. \item $S= \s_\ga^{-1}(A^*)$. \item $A^* = \langle \pi_\ga (S)
, \pi_\ga (S)^{-1} \rangle$, i.e. the group of units of the ring
$A$ is generated by the sets $\pi_\ga (S)$ and $\pi_\ga^{-1}(S):=
\{ \pi_\ga (s)^{-1} \, | \, s\in S\}$. \item $A^* = \{ \pi_\ga
(s)^{-1}\pi_\ga ( t) \, |\, s, t\in S\}$. \item $Q_l(A) = A$ and
$\Ass_l(A) = \{ 0\}$.     In particular, if $T\in \Den_l(A, 0)$
then  $T\subseteq A^*$.
\end{enumerate}
\end{theorem}


Let $\maxLoc_l(R)$ be the set of maximal elements of the poset
$(\Loc_l(R), \ra )$ where $A\ra B$ for $A,B\in \Loc_l(R)$ means that there exist $S,T\in \Den_l(R)$ such that $S\subseteq T$, $A= S^{-1}R$  and $B= T^{-1}R$ (then there exists a natural ring homomorphism $A\ra B$, $s^{-1}r\mapsto s^{-1}r$). Then (see \cite{larglquot}),
\begin{equation}\label{mADen1}\maxLoc_l(R) = \{ S^{-1}R \, | \, S\in \maxDen_l(R) \}= \{ Q_l(R/
\ga ) \, | \, \ga \in \assmaxDen_l(R)\}.
\end{equation}


{\bf The maximal left quotient rings of a finite direct product of rings}.
\begin{theorem}\label{c26Dec12}
 \cite{Bav-Crit-S-Simp-lQuot} Let  $R=\prod_{i=1}^n R_i$ be the direct product of rings $R_i$. Then
for each $i=1, \ldots , n$, the map
\begin{equation}\label{aab1}
\maxDen_l(R_i) \ra \maxDen_l(R), \;\; S_i\mapsto R_1\times\cdots \times S_i\times\cdots \times R_n,
\end{equation}
is an injection. Moreover, $\maxDen_l(R)=\coprod_{i=1}^n \maxDen_l(R_i)$ in the sense of (\ref{aab1}), i.e.
$$ \maxDen_l(R)=\{ S_i\, | \, S_i\in \maxDen_l(R_i), \; i=1, \ldots , n\},$$
$S_i^{-1}R\simeq S_i^{-1}R_i$, $\ass_R(S_i)= R_1\times \cdots \times \ass_{R_i}(S_i)\times\cdots \times R_n$. The core of the left denominator set $S_i$ in $R$ coincides with the core of the left denominator set $S_i$ in $R_i$, i.e.
$$(R_1\times\cdots \times S_i\times\cdots \times R_n)_c=0\times\cdots \times S_{i,c}\times\cdots \times 0.$$
\end{theorem}

\begin{corollary}\label{b16Feb14}
Let $R= \prod_{i=1}^nR_i$ be the direct product of rings $R_i$. Then $\gll_R=\prod_{i=1}^n\gll_{R_i}$.
\end{corollary}

{\bf A bijection between $\maxDen_l(R)$ and $\maxDen_l(Q_l(R))$}.
\begin{proposition}\label{A8Dec12}
  \cite{Bav-Crit-S-Simp-lQuot} Let $R$ be a ring, $S_l$ be the  largest regular left Ore set of the ring $R$, $Q_l:= S_l^{-1}R$ be the largest left quotient ring of the ring $R$, and $\CC$ be the set of regular elements of the ring $R$. Then
\begin{enumerate}
\item $S_l\subseteq S$ for all $S\in \maxDen_l(R)$. In particular,
$\CC\subseteq S$ for all $S\in  \maxDen_l(R)$ provided $\CC$ is a
left Ore set. \item Either $\maxDen_l(R) = \{ \CC \}$ or,
otherwise, $\CC\not\in\maxDen_l(R)$. \item The map $$
\maxDen_l(R)\ra \maxDen_l(Q_l), \;\; S\mapsto SQ_l^*=\{ c^{-1}s\,
| \, c\in S_l, s\in S\},
$$ is a bijection with the inverse $\CT \mapsto \s^{-1} (\CT )$
where $\s : R\ra Q_l$, $r\mapsto \frac{r}{1}$, and $SQ_l^*$ is the
sub-semigroup of $(Q_l, \cdot )$ generated by the set  $S$ and the
group $Q_l^*$ of units of the ring $Q_l$, and $S^{-1}R= (SQ_l^*)^{-1}Q_l$.
    \item  If $\CC$ is a left Ore set then the map $$ \maxDen_l(R)\ra \maxDen_l(Q:= \CC^{-1}R), \;\; S\mapsto SQ^*=\{ c^{-1}s\,
| \, c\in \CC, s\in S\}, $$ is a bijection with the inverse $\CT
\mapsto \s^{-1} (\CT )$ where $\s : R\ra Q$, $r\mapsto
\frac{r}{1}$, and $SQ^*$ is the sub-semigroup of $(Q, \cdot )$
generated by the set  $S$ and the group $Q^*$ of units of the ring
$Q$, and $S^{-1}R= (SQ^*)^{-1}Q$.
\end{enumerate}
\end{proposition}


\section{Idempotent left denominator sets}\label{IDEMLDS}

Theorem \ref{A9Feb13} states that every left localization of a left Artinian ring is an idempotent localization.
In this section, several results on idempotent left denominator sets are given that are used in proofs of the subsequent sections.

 Let $R$ be a ring and  $e\in R$ be a {\em nonzero} idempotent. Then $1=e_1+e_2$ is the sum of orthogonal idempotents where $e_1=e$ and $e_2=1-e_1$. The ring $R$ can be seen as the {\em matrix ring associated with the idempotent} $e$,
\begin{equation}\label{MRawe}
R=\bigoplus_{i,j=1}^2R_{ij}= \begin{pmatrix}
 R_{11}& R_{12}\\ R_{21} & R_{22}
 \end{pmatrix}\;\; {\rm where}\;\; R_{ij}:=e_iRe_j.
\end{equation}
For an element $r\in R$, let $r\cdot :R\ra R$, $x\mapsto rx$, and $\cdot r :R\ra R$, $x\mapsto xr$.  The next proposition is a criterion for an {\em idempotent multiplicative  set} $S_e=\{ 1,e\}$ to be a left denominator set of the ring $R$.

 \begin{proposition}\label{c24Dec12}
Let $e$ be a nonzero idempotent of  a ring $R$. We keep the notation as above. Then $S_e=\{ 1,e\}\in \Den_l(R)$ iff $R_{12}=0$. In this case, $\ass (S_e)=\ker (e\cdot ) = (1-e)R$ and $S_e^{-1}R\simeq R/\ass (S_e)\simeq R/ (1-e)R\simeq R_{11}$.
\end{proposition}

{\it Proof}. $(\Rightarrow )$ If $S_e\in \Den_l(R)$ then $\ga := \ass (S_e)=\ker (e_1\cdot ) = e_2R$. Since $R_{12}e_1=0$, we must have $R_{12}\subseteq \ga = e_2R= R_{21}+R_{22}$, hence $R_{12}=0$. Clearly, $S_e^{-1}R\simeq R/\ga \simeq R_{11}$ since $e_2\in \ga$ and $1\equiv e_1\mod \ga$.

$(\Leftarrow )$ Suppose that $R_{12}=0$, i.e. $R=\begin{pmatrix}
 R_{11}& 0\\ R_{21} & R_{22}
 \end{pmatrix}$. Then $S_e\in \Ore_l(R)$ since for any element $r= \begin{pmatrix}
 r_{11}& 0\\ r_{21} & r_{22}
 \end{pmatrix}\in R$, $er= \begin{pmatrix}
 r_{11}& 0\\ 0 & 0
 \end{pmatrix}= \begin{pmatrix}
 r_{11}& 0\\ 0 & 0
 \end{pmatrix} e$.  The inclusion  $\ker (\cdot e)= R_{22}\subseteq \ga$ implies that $S_e\in \Den_l(R)$.  $\Box $

$\noindent $

 Proposition \ref{c24Dec12} means that the idempotent multiplicative set $S_e$ is a left denominator set of $R$ iff the ring $R$  is {\em left triangular} (as the matrix ring associated with the idempotent $e$), i.e.  $$R=\begin{pmatrix}
 R_{11}& 0\\ R_{21} & R_{22}
 \end{pmatrix}.$$

This fact is the most vivid demonstration of the fact that very often a left denominator set fails to be a right denominator set. Recall that $\Den (R)$ is the set of (left and right) denominator sets of a ring $R$. The next corollary demonstrates that the condition being a {\em  left and right} denominator set  is a strong one.

\begin{corollary}\label{d24Dec12}
Let $e$ be a nonzero idempotent of a ring $R$. Then $S_e=\{ 1,e\}\in \Den (R)$ iff $R_{12}=0$ and $R_{21}=0$ iff $e$ is a central idempotent. In this case, $\ass (S_e)=\ker (e\cdot ) = (1-e)R$ and $S_e^{-1}R\simeq R/ (1-e)R\simeq R_{11}$.
\end{corollary}

{\it Proof}. The first `iff' is due to Proposition \ref{c24Dec12}. The second `iff' is obvious.  $\Box $

$\noindent $

Let $\Aut (R)$ be the group of automorphisms of the ring $R$ and $\Inn (R):=\{ \o_u\, | \, u\in R^*\}$ be the group of inner automorphisms of the ring $R$ where $\o_u (r) := uru^{-1}$ for $r\in R$. The group $\Inn (R)$ is a normal subgroup of $\Aut (R)$ (since for $\s \in \Aut (R)$ and $\o_u\in \Inn (R)$, $\s\o_u\s^{-1} = \o_{\s (u)}$).

Let $e\in R$ be an idempotent. If $S_e=\{ 1,e\}\in \Den_l(R)$ then the set $S_e$ is called an {\em idempotent left denominator set} of the ring $R$ and the idempotent $e$ is called a {\em left denominator  idempotent} of $R$. Let $\IDen_l(R)$ be the set of all the idempotent left denominator sets of the ring $R$ and let $\CI_l = \CI_l(R)$ be the set of all left  denominator  idempotents of the ring $R$. By Proposition \ref{c24Dec12},
\begin{equation}\label{IlR=d}
\CI_l(R)=\{ e\in R\, | \, e^2=e, \, eR(1-e)=0\}.
\end{equation}
The map
\begin{equation}\label{CIID}
\CI_l(R)\ra \IDen_l(R), \;\;\ e\mapsto S_e=\{ 1,e\},
\end{equation}
is a bijection.  The groups $\Aut (R)$ and $\Inn (R)$ act in the obvious way on the sets $\CI_l(R)$ and $\IDen_l(R)$. A ring $R$ is called a {\em local ring} if the factor ring $R/ \rad (R)$ is a division ring.

\begin{corollary}\label{b14Feb13}
Let $R$ be a ring and $e\in R\backslash \{ 0,1\}$ be an idempotent.
\begin{enumerate}
\item The following statements are equivalent.
\begin{enumerate}
\item The idempotents $e$ and $1-e$ are left  denominator  idempotents.
\item The idempotent $e$ is a central idempotent.
\item The idempotents $e$ and $1-e$ are right   denominator  idempotents.
\end{enumerate}
\item The following statements are equivalent.
\begin{enumerate}
\item All the idempotents of $R$ are left  denominator  idempotents.
\item All the idempotents of $R$ are central idempotents.
\item All the idempotents of $R$  are right   denominator  idempotents.
\end{enumerate}

\item Let $R$ be a left Artinian ring. Then the  following statements are equivalent.

\begin{enumerate}
\item All the idempotents of $R$ are left  denominator  idempotents.
\item The ring $R$ is a direct product of finitely many  local left Artinian rings.
\item All the idempotents of $R$ are are right   denominator  idempotents.
\end{enumerate}
\end{enumerate}
\end{corollary}

{\it Proof}. 1. Statement 1 follows from Proposition \ref{c24Dec12}.

2. Statement 2 follows from statement 1.

3. Statement 3 follows from statement 2 and the fact that 1 is the only nonzero idempotent of a local left Artinian ring (in such a ring, every nonzero idempotent $e$ is a primitive one, hence is conjugate to 1, i.e. $e=1$). $\Box $

$\noindent $

A set $(X, \geq )$ is called a {\em pre-ordered set} if

(i) $x\geq x$,

(ii) $x\geq y$ and $y\geq z$ implies $x\geq z$.

In general, the conditions $x\geq y$ and $y\geq x$ do not imply $x=y$. If this property holds the pre-ordered set $X$ is called a {\em partially ordered set}, a {\em poset}, for short. The set $\CI_l$ is a pre-ordered  set $(\CI_l, \geq )$ where
$$e_1\geq e_2\;\; {\rm  iff}\;\;  e_2e_1=e_2\;\; {\rm  iff}\;\; (1-e_2)(1-e_1)=1-e_1\;\; {\rm  iff}\;\; (1-e_1)R\subseteq (1-e_2)R.$$

 The last `iff'  follows from the fact that the inclusion $(1-e_1)R\subseteq (1-e_2)R$ implies the inclusion $e_2(1-e_1)R\subseteq e_2 (1-e_2)R=0$, and so $e_2e_1=e_2$.   Via the bijection (\ref{CIID}), the set $\IDen_l(R)$ is a pre-ordered set $(\IDen_l(R), \geq )$ where $S_{e_1}\geq S_{e_2}$ iff $e_1\geq e_2$.

$\noindent $

{\it Definition}.  We say that a left localization $A=S^{-1}R$  of a ring $R$ is an {\em idempotent left localization} of $R$ if there is an idempotent $e\in R$ such that $S_e= \{ 1,e\}\in \Den_l(R)$ and   $S_e^{-1}R= A$ in $\Loc_l(R)$, i.e. the map $A= S^{-1}R\ra S_e^{-1}R$, $ s^{-1}r\mapsto s^{-1}r$, is an isomorphism, i.e. $\ass (S) = \ass (S_e)$ and $s+\ass (S_e)\in ( R/ \ass (S_e))^*$ for all $s\in S$. Let $\ILoc_l(R):=\{ S_e^{-1}R\, | \, e\in \CI_l(R)\}$, the set of all the idempotent left localizations of $R$.

$\noindent $

{\it Remark}. In general, even for Artinian rings, it is not true that the condition $S^{-1}R\simeq S_e^{-1}R$ or even $ S^{-1}R = S_e^{-1}R$ in $\Loc_l(R)$ for some $S\in \Den_l(R)$ and $S_e\in \IDen_l(R)$ implies that  the set $S$  contains an idempotent element distinct from 1.

Example: Let $R=\begin{pmatrix}
\Q & 0 \\ \Q & \Q
 \end{pmatrix}$, $e=E_{11}$, $\ga = (1-e)R$, $s=\begin{pmatrix}
 2& 0\\ 1 & 0
 \end{pmatrix}$, $S=\{1, s^i=\begin{pmatrix}
2^i & 0\\ 2^{i-1} & 0
 \end{pmatrix}\, | \, i\geq 1\}$. Then $S_e=\{1,e\}\in \Den_l(R, \ga )$ (by Proposition \ref{c24Dec12}) and $S\in \Den_l(R, \ga )$ (by Corollary \ref{b14Feb13}),  $S_e^{-1}R\simeq R/ \ga \simeq S^{-1}R$ (by Proposition \ref{c24Dec12} and Corollary \ref{b14Feb13}). Clearly, 1 is the only idempotent of the set $S$ and $S_e^{-1}R = S^{-1}R$ in $\Loc_l(R)$. For all $i\geq 1$, $s^i\ga =0$ and $\ga s^i\neq 0$.

\begin{lemma}\label{a11Feb13}
Let $e_1,e_2\in \CI_l(R)$. The following statements are equivalent.
\begin{enumerate}
\item $e_1\geq e_2$.
\item $\frac{e_1}{1}$ is a unit in the ring $S_{e_2}^{-1}R$;  equivalently, $\frac{e_1}{1}=1\in S_{e_2}^{-1}R$.
\item The map $S_{e_1}^{-1}R\ra S_{e_2}^{-1}R$, $e_1^ir\ra e_1^ir$, where $i=0,-1$, is well-defined.
\item $\ass (S_{e_1})\subseteq \ass (S_{e_2})$.
\end{enumerate}
\end{lemma}

{\it Proof}. $(1\Rightarrow 2)$ If $e_1\geq e_2$, i.e. $e_2e_1=e_2$ then $\frac{e_2}{1}\frac{e_1}{1}=\frac{e_2}{1}$ in $S_{e_2}^{-1}R$, and so $\frac{e_1}{1}=1=\frac{e_2}{1}$ since $\frac{e_2}{1}$ is the identity element of the ring $S_{e_2}^{-1}R$.

$(2\Rightarrow 3)$ If $\frac{e_1}{1}$ is a unit of the ring $S_{e_2}^{-1}R$ then the map $e_1^ir\ra e_1^ir$ exists by the universal property of left localizations.

$(3\Rightarrow 1)$ If the map in statement 3 is well-defined then $\frac{e_1}{1}=1=\frac{e_2}{1}$ in $S_{e_2}^{-1}R$ and so $e_1-e_2\in \ass (S_{e_2})= \ker_R(e_2\cdot )$, i.e. $e_2(e_1-e_2)=0$. This means that $e_1\geq e_2$.

$(1\Rightarrow 4)$  The equality   $e_2e_1=e_2$  implies that        $\ass (S_{e_1})=\ker (e_1\cdot ) \subseteq \ker (e_2\cdot )=\ass (S_{e_2})$.

$(4\Rightarrow 2)$ Notice that $S_{e_1}^{-1}R=R/\ass (S_{e_1})$, $ S_{e_2}^{-1}R=R/\ass (S_{e_2})$ and $e_1\in 1+\ass (S_{e_1})$. If $\ass (S_{e_1})\subseteq \ass (S_{e_2})$ then $S_{e_2}^{-1}R\ni \frac{e_1}{1}= 1+\ass (S_{e_2})$ is a unit.  $\Box $

$\noindent $

Let us define an equivalence relation $\sim$ on $\CI_l(R)$ by the rule $e_1\sim e_2$ iff $e_1\geq e_2$ and $e_2\geq e_1$. Let $[e]:=\{ f\in \CI_l(R)\, | \, f\sim e\}$ be the equivalence class of $e\in \CI_l(R)$. Then set $\CI_l(R)/\sim := \{ [e]\, | \, e\in \CI_l(R)\}$ of equivalence classes is a poset $(\CI_l(R)/\sim , \geq )$ where $[e_1]\geq [e_2]$ if $e_1\geq e_2$.

\begin{lemma}\label{b11Feb13}
Let $e_1,e_2\in \CI_l(R)$. The following statements are equivalent.
\begin{enumerate}
\item $e_1\sim e_2$.
\item $e_1\equiv 1\mod \ass (S_{e_2})$ and $e_2\equiv 1\mod \ass (S_{e_1})$ (i.e. $S_{e_1}^{-1}R = S_{e_2}^{-1}R$ in $\Loc_l(R)$).
\item The map   $S_{e_1}^{-1}R\ra  S_{e_2}^{-1}R$, $e_1^ir\mapsto e_1^ir$, where $i=0,-1$, is an isomorphism.
\item $\ass (S_{e_1})=\ass (S_{e_2})$.
\end{enumerate}
\end{lemma}

{\it Proof}. The lemma is an easy corollary of Proposition \ref{c24Dec12} and Lemma \ref{a11Feb13}. $\Box$

Notice that every ideal is invariant under the inner automorphisms. By Lemma \ref{b11Feb13}, the group of inner automorphisms $\Inn (R)$ of the ring $R$ acts on the set $\CI_l(R)/\sim$ by the rule: for any $u\in R^*$ and $[e]\in \CI_l(R)/\sim$, $\;\; u[e]u^{-1}:= [ueu^{-1}]$.

$\noindent $

{\it Example}. Let $R= \begin{pmatrix}
R_{11} & 0 \\ R_{21} & R_{22}
 \end{pmatrix}$
  be any triangular ring where $R_{11}$ and $R_{22}$ are arbitrary rings and $R_{21}$ be an arbitrary $(R_{22},R_{11})$-bimodule. Let $e_1=E_{11}$. By Proposition \ref{c24Dec12}, $S_{e_1}=\{ 1,e_1\}\in \Den_l(R, \ga )$ where $\ga = E_{22}R$. Clearly, $R^*= \begin{pmatrix}
R_{11}^* & 0 \\ R_{21} & R_{22}^*
 \end{pmatrix}$ and $\Inn (R)\cdot e_1= \begin{pmatrix}
1 & 0 \\ R_{21} & 0
 \end{pmatrix}$ since for all units $u= \begin{pmatrix}
\alpha  & 0 \\ \beta & \gamma
 \end{pmatrix}\in R^*$,
 $$ ue_1u^{-1}=\begin{pmatrix}
1 & 0 \\ \beta\alpha^{-1} & 0
 \end{pmatrix}.$$
 Therefore, for any element $a\in R_{21}$,
$e_a:= \begin{pmatrix}
1 & 0 \\ a & 0
 \end{pmatrix}$  is an idempotent of $R$, $S_{e_a}\in \Den_l(R, \ga )$ and $S_{e_a}^{-1}R\simeq R/ \ga =S^{-1}_{e_1}R$.


\section{Left localizations of left Artinian rings}\label{LLLAR}

Throughout this section, $R$ {\em is a left Artinian ring if it is  not stated otherwise.} The aim of this section is to  prove, for a left Artinian ring $R$, that every left localization of $R$ is an idempotent left localization (Theorem \ref{A9Feb13}), there are only finitely many left localization rings of $R$ and to give a classification of all of them (Theorem \ref{A11Feb13}), to classify the maximal left denominator sets of $R$ (Theorem \ref{B11Feb13}), to give an explicit description of the left localization radical $\gll_R$ of $R$ (Theorem \ref{D12Feb13}). The ideals in $\Ass_l(R)$ have many interesting/unexpected properties (Corollary \ref{A13Feb13}).

An element $r$ of a ring $R$ is called a {\em left regular} if the map $\cdot r:R\ra R$ is an injection, i.e.   $xr=0$ implies $x=0$. Recall that every left Artinian ring is left Noetherian.

\begin{lemma}\label{a9Feb13}
Let $R$ be a left Artinian ring.
\begin{enumerate}
\item Every left regular element of the ring $R$ is a unit.
\item Let $S\in \Den_l(R, \ga )$. Then $S^{-1}R\simeq R/ \ga$ (an $R$-isomorphism).
\item Let $s\in R$ and $\cdot s:R\ra R$, $r\mapsto rs$. Then $R=Rs\oplus \ker (\cdot s)$ iff $Rs\cap \ker (\cdot s)=0$.
    \item Let $S\in \Den_l(R, \ga )$ and $T\in \Den_l(R, \gb )$. Then the rings $S^{-1}R$ and $T^{-1}R$ are $R$-isomorphic iff $\ga = \gb$.
\end{enumerate}
\end{lemma}

{\it Proof}. 1. Trivial.

2. Statement 2 follows from statement 1. Let $\pi : R\ra R/\ga$, $r\mapsto \br = r+\ga$. Then $\pi (S)\in \Den_l(R/\ga , 0)$ and $S^{-1}R\simeq \pi (S)^{-1}(R/\ga )$. By statement 1, the set of regular elements $\pi (S)$ of the left Artinian ring $R/\ga $ consists of units, and so $\pi (S)^{-1}(R/ \ga )\simeq R/\ga$.

3. $(\Rightarrow )$ Trivial.

$(\Leftarrow )$ if $Rs\cap \ker (\cdot s )=0$ then $Rs\oplus \ker (\cdot s)\subseteq R$. It follows from the short exact sequence of $R$-modules $0\ra \ker (\cdot s)\ra R\stackrel{\cdot s}{\ra } Rs\ra 0$ that the left $R$-modules $Rs\oplus \ker (\cdot s)$ and $ R$ have the same length. Then, $Rs\oplus \ker (\cdot s)=R$.

4. Statement 4 follows from statement 2. $\Box $

$\noindent $

Statement 4 does not hold for non-Artinian rings, eg, $R=\Z$, $S_1=\{ 1\}$ and $S_2= \{ 2^i\, | \, i\in \N\}$. Statement 4 does not hold  if the condition  `$R$-isomorphic' is replaced by `isomorphic', eg, if $R= K\times K$, $1= e_1+e_2$, $K$ is a field then $S_{e_1}^{-1}R \simeq K \simeq S_{e_2}^{-1}R$ but $\ass (S_{e_1})=(1-e_1)R= \{ 0\} \times K \neq K\times \{ 0\} = (1-e_2)R=\ass (S_{e_1})$.

{\bf The idempotent $e(s)$ associated with $s\in R$}.
Suppose that $R$ is a left Artinian ring. Then for each element $s\in R$ there is the least natural number $n=n(s)$ such that $Rs^n= Rs^i$ for all $i\geq n$. By Lemma \ref{a9Feb13}.(3), the number $n=n(s)$ is the least natural number $n$ such that $Rs^i\cap \ker (\cdot s^i)=0$ for all $i\geq n$ or equivalently $R=Rs^i\oplus \ker (\cdot s^i)$ for all $i\geq n$ (equivalently, $\ker (\cdot s^i) =\ker (\cdot s^n)$ for all $i\geq n$). So, for all $i\geq n$, $Rs^n=Rs^i$ and $\ker(\cdot s^n) = \ker (\cdot s^i)$. For each  $i\geq n$, let
\begin{equation}\label{1=ese}
1=e[s^i]+e'[s^i]
\end{equation}
be the sum of orthogonal idempotents that corresponds to the decomposition $R=Rs^i\oplus \ker (\cdot s^i)$. Then $e[s^n]=e[s^i]$ and $e'[s^n]=e'[s^i]$ for all $i\geq n$. These common values are denoted by $e(s)$ and $e'(s)$, respectively. The idempotent $e(s)$ of the ring $R$ is called the {\em idempotent associated with the element} $s$. The element $s$ is a nilpotent iff $e(s)=0$. The element $s$ is a unit  iff $e(s)=1$. The element $s$ is neither a nilpotent element nor a unit iff $e(s)\neq 0,1$.

{\bf Orthogonal idempotents}.
The ring  $R$ is a left Artinian ring. Its radical $\rad (R)$ is a nilpotent ideal, and so it coincides with the prime radical $\gn_R$ of the ring $R$, and
\begin{equation}\label{bRPMR}
\bR :=R/\rad (R) \simeq \prod_{i=1}^s M_{n_i}(D_i)
\end{equation}
where $\bR_i:= M_{n_i}(D_i)$ is the ring of $n_i\times n_i$ matrices with entries from a division ring $D_i$. Let $\{ \bE_{pq}(i)\, | \, p,q=1, \ldots  , n_i\}$ be the matrix units of the ring $\bR_i$. Since the radical $\rad (R)$ is a nil ideal (even a nilpotent ideal) of $R$, the decomposition of 1 in $\bR$ as a sum of primitive orthogonal idempotents
\begin{equation}\label{1=SPI}
1=\sum_{i=1}^s\sum_{j=1}^{n_i} \bE_{jj}(i)
\end{equation}
can be lifted to a decomposition of $1\in R$ as a sum of primitive orthogonal idempotents
\begin{equation}\label{1=SPI1}
1=\sum_{i=1}^s\sum_{j=1}^{n_i} E_{jj}(i)
\end{equation}
and any such a lift is unique up to conjugation (i.e. up to inner automorphism) and order of idempotents (i.e. permutation of the idempotents), (Proposition 18.23.5, \cite{Faith-Book-II}).    The sum $\overline{1}_i:=\sum_{j=1}^{n_i}\bE_{jj}(i)$ is the identity of the ring $\bR_i$ and let $1_i:=\sum_{j=1}^{n_i}E_{jj}(i)$. Then $1=\sum_{i=1}^s\overline{1}_i$ is the sum of central  orthogonal idempotents of the ring $\bR$, and
\begin{equation}\label{1=S1is}
1=\sum_{i=1}^s1_i
\end{equation}
is the sum of orthogonal idempotents of the ring $R$ (in general, not necessarily central). For each non-empty subset $I$ of the set $\{ 1, \ldots , s\}$, let
\begin{equation}\label{1I=1i}
e_I:=1_I:=\sum_{i\in I}1_i.
\end{equation}
Since the idempotents $\overline{1}_i$  are central, the sum (\ref{1=S1is}) is {\em unique up to inner automorphism} of $R$, i.e. if $1=\sum_{i=1}^s 1_i'$  is another lifting with $\overline{1}_i'= \overline{1}_i$ then $1_i'= u1_iu^{-1}$ for all $i$ and some unit $u\in R^*$.
 Let $\ga$ be an ideal of the ring $R$ such that $\ga \not\subseteq \rad (R)$. Then $(\ga +\rad (R))/ \rad (R)$ is a nonzero ideal of the semisimple ring $\bR$, and so
\begin{equation}\label{arMn}
(\ga +\rad (R))/ \rad (R)= \prod_{i\in  I(\ga )}M_{n_i}(D_i)=\bR \overline{1}_{I(\ga )}
\end{equation}
for a unique non-empty subset $I(\ga )$ of the set $\{ 1, \ldots , s\}$ where $\overline{1}_{I(\ga )}:=\sum_{i\in I(\ga )}\overline{1}_i$ is the central idempotent of the ring $\bR$. The set $I(\ga )$ is a proper subset of $\{ 1, \ldots , s\}$ iff $\ga \neq R$.

{\bf Every left localization of a left Artinian ring is an idempotent left localization}.  Proposition \ref{G12Feb13} and Corollary \ref{d12Feb13} are about lifting (in many different ways) left denominator sets of a factor ring of a ring to the ring (under certain conditions).

\begin{proposition}\label{G12Feb13}
Let $R$ be an arbitrary ring, $\ga$ be its ideal, $\bR := R/\ga $ and $S$ be a multiplicative set of the ring $R$ such that for each element $a\in \ga$ there is an element $s\in S$ such that $sa=0$.   Let $\pi : R\ra \bR$, $r\mapsto \br = r+\ga$, and $\bS := \pi (S)$.
\begin{enumerate}
\item If $\bS\in \Ore_l(\bR , \overline{\gb})$ then $S\in \Ore_l(R, \gb )$ where $\gb = \pi^{-1} ( \overline{\gb})$.
\item If $\bS\in \Den_l(\bR , \overline{\gb})$ then $S\in \Den_l(R, \gb )$ where $\gb = \pi^{-1} ( \overline{\gb})$ and $S^{-1}R \simeq \bS^{-1}\bR$.
\end{enumerate}
\end{proposition}

{\it Proof}. 1. (i) $S\in \Ore_l(R )$: For given  elements $s\in S$ and $r\in R$, we have to find elements $s'\in S$ and $r'\in R$ such that $s'r=r's$. Since $\bS\in \Ore_l(\bR )$, $\bs_1\br = \br_1\bs$ for some elements $s_1\in S$ and $r_1\in R$. Then $s_1r-r_1s\in \ga$, and so $s_2(s_1r-r_1s)=0$ for some element $s_2\in S$. It suffices to take $s' = s_2s_1$ and $r=s_2r_1$.

(ii) $\ass (S) = \gb$: Let $sr=0$ for some elements $s\in S$ and $r\in R$. Then $\bs \br =0$ in $\bR$, and so $\br\in \overline{\gb}$, hence $r\in \gb$, i.e. $\ass (S)\subseteq \gb$.

Let $b\in \gb$. Then $\bs \ob =0$ for some element $s\in S$. Then $sb\in \ga$, and so $s_1sb=0$ for some element $s_1\in S$. Hence $\gb \subseteq \ass (S)$.

2. (i) $S\in \Den_l(R, \gb )$:  In view of statement 1, it suffices to show that $rs=0$ for some elements $s\in S$ and $r\in R$ implies $s'r=0$ for some element $s'\in S$. We have the equality $\br \bs =0$ in the ring $\bR$. Then $\bs_1\br =0$  for some element $s_1\in S$ since $\bS\in \Den_l(R, \overline{\gb} )$. Then $s_1r\in \ga$, hence $s_2s_1r=0$ for some element $s_2\in S$. It suffices to take $s'=s_2s_1$.

(ii) $S^{-1}R\simeq \bS^{-1}\bR$: By the universal property of left localizations, the map $S^{-1}R\ra \bS^{-1}\bR$, $s^{-1}r\mapsto \bs^{-1}\br$, is a ring homomorphism which is obviously an epimorphism. Suppose that an element $s^{-1}r\in S^{-1}R$ belongs to the kernel of the epimorphism. Then $\br =0$ in $\bS^{-1}\bR$, and so $\bs_1\br =0$ in $\bR$ for some element $s_1\in S$. This means that $s_1r\in \ga$, and so $s_2s_1r=0$ for some element $s_2\in S$. Therefore, $\frac{r}{1}=0$ and $s^{-1}r=0$. The epimorphism is an isomorphism.  $\Box $


\begin{corollary}\label{d12Feb13}
Let $R$ be a  ring and $S\in \Den_l(R, \ga )$ (respectively,  $S\in \Ore_l(R, \ga )$). Then $S+\ga \in \Den_l(R, \ga )$ (respectively, $S+\ga \in \Ore_l(R, \ga )$) and $(S+\ga )^{-1}R\simeq S^{-1}R$.
\end{corollary}

{\it Proof}. Since $S\cap \ga =\emptyset$, the set $S+\ga$ is a multiplicative set.  We keep the  notation of Proposition \ref{G12Feb13}. Since $S\in \Den_l(R, \ga )$ (respectively,  $S\in\Ore_l(R, \ga )$), we have that $\overline{S+\ga}=\bS \in \Den_l(\bR , 0)$ (respectively,  $\bS\in\Ore_l(\bR, 0 )$). By Proposition \ref{G12Feb13}, $S+\ga \in \Den_l(R, \ga )$ (respectively, $S+\ga \in \Ore_l(R, \ga )$) and $(S+\ga )^{-1}R\simeq \bS^{-1}\bR\simeq S^{-1}R$.  $\Box $


\begin{corollary}\label{b14Feb13}
Let $R$ be a ring, $\ga$ be an ideal of $R$ such that $\ga\neq R$, an element $s\in R$ be such that the element $s+\ga$ is a unit of the ring $R/\ga$ and $\bigcup_{i\geq 1}\ker_R(s^i\cdot )=\ga$. Then $S_s:=\{ s^i\, | \, i\in \N\}\in \Den_l(R, \ga )$ and $S_s^{-1}R \simeq R/ \ga$.
\end{corollary}

{\it Proof}. We keep the  notation of Proposition \ref{G12Feb13}.
 The set $\bS_s:=\{ \bs^i\, | \, i\in \N\}$ (where $\bs = s+\ga$) consists of units of the ring $\bR = R/ \ga$, and so $\bS_s\in \Den_l(\bR , 0)$. By Proposition \ref{G12Feb13}, $S_s\in \Den_l(R, \ga )$. Now, it is obvious that  $S_s^{-1}R \simeq R/ \ga$ (since $\bS_s\subseteq \bR^*$).  $\Box $

$\noindent $

Let $S$ be a nonempty subset of a ring $R$. The set $\ker_l(S)=\{ \ker (\cdot s)\, | \, s\in S\}$ of left ideals of the ring $R$ is a poset with respect to $\subseteq$. Let $\maxker_l(S)$ be the set of maximal elements of the poset $\ker_l(S)$.  The set $\maxker_l(S)$ is a non-empty set provided the ring $R$ satisfies the a.c.c. for  left annihilators. Similarly, the set $\ker_r(S)=\{ \ker ( s\cdot )\, | \, s\in S\}$ of right  ideals of the ring $R$ is a poset with respect to $\subseteq$. Let $\maxker_r(S)$ be the set of maximal elements of the poset $\ker_r(S)$.  The set $\maxker_r(R)$ is a non-empty set provided the ring $R$ satisfies the a.c.c. for  right annihilators.

The first statement of the following theorem shows that every left localization of a left Artinian ring is an idempotent left localization.
\begin{theorem}\label{A9Feb13}
Let $R$ be a left Artinian ring and $S\in \Den_l(R, \ga )$. Then
\begin{enumerate}
\item There exists a nonzero idempotent $e\in R$ such that $S_e:=\{ 1,e\}\in \Den_l(R, \ga )$ and the rings $S^{-1}R$ and $ S_e^{-1}R$ are $R$-isomorphic.
\item
\begin{enumerate}
\item If $\ga =0$ then $e=1$.
\item If $\ga \neq 0$ then $\ga = (1-e)R=1_{I(\ga )}R\not\subseteq \rad (R)$ and the idempotent $e$ is conjugate to $1_{CI (\ga )}=1-1_{I(\ga )}$ where the set $I(\ga )$ is defined in (\ref{arMn}).
    \item 
\begin{equation}\label{R=Rij1}
R=\begin{pmatrix}
 R_{11}& 0\\ R_{21} & R_{22}
 \end{pmatrix}, \;\;\; \ga =\begin{pmatrix}
 0&0 \\ R_{21} & R_{22}
 \end{pmatrix} \;\; {\rm and}\;\; S\subseteq \begin{pmatrix}
 R_{11}^*& 0\\ R_{21} & R_{22}
 \end{pmatrix},
\end{equation}
where $R_{ij}=e_iRe_j$, $e_1=e$ and $e_2=1-e_1$.
\end{enumerate}
\end{enumerate}
\end{theorem}

{\it Proof}. 1. If $\ga =0$ then $S\subseteq R^*$ (Lemma \ref{a9Feb13}.(1)) and so $S^{-1}R = R$. It suffices to take $e=1$.

We can assume that $\ga \neq 0$. Then necessarily $\ga \neq R$, i.e. $\ga$ is  a proper ideal of the ring $R$ and $S\not\subseteq R^*$.
By Corollary \ref{d12Feb13}, $S+\ga \in \Den_l(R, \ga )$ and $(S+\ga )^{-1} R\simeq S^{-1}R$ (an $R$-isomorphism). Without loss of generality we may assume that $S+\ga \subseteq S$ (replacing $S$ by $S+\ga$, if necessary).

The ring $R$ is a left Artinian ring. Therefore, we can fix an element $s\in S$ such that $\ker ( \cdot s)\in \maxker_l(S)$, and so $\ker (\cdot s)\neq 0$ (since $S\not\subseteq R^*$ and by Lemma \ref{a9Feb13}.(1)). Then necessarily $Rs\cap \ker (\cdot s)=0$, and so
\begin{equation}\label{R=Rspk}
R=Rs\oplus \ker (\cdot s)
\end{equation}
is the direct sum of left ideals of the ring $R$, by Lemma \ref{a9Feb13}.(3). Let $1=e_1+e_2$ be the corresponding decomposition of 1 as a sum of orthogonal idempotents. Then
$$ R=\bigoplus_{i,j=1}^2R_{ij}= \begin{pmatrix}
 R_{11}& R_{12}\\ R_{21} & R_{22}
 \end{pmatrix}\;\; {\rm where}\;\; R_{ij}:=e_iRe_j, $$
$Rs=Re_1=\begin{pmatrix}
 R_{11}& 0\\ R_{21} & 0
 \end{pmatrix}$ and $ \ker (\cdot s)= Re_2= \begin{pmatrix}
 0& R_{12}\\ 0 & R_{22}
 \end{pmatrix}.$ Since $s\in Rs$, the element $s$ has the form
\begin{equation}\label{s=00}
s=\begin{pmatrix}
 s_{11}& 0\\ s_{21} & 0
 \end{pmatrix},
\end{equation}
 the ring $R_{11}$ is a left Artinian ring and the map $\cdot s: R_{11}\ra R_{11}$ is an injection, by (\ref{R=Rspk}). By Lemma \ref{a9Feb13}, $s_{11}$ is a unit of the ring $R_{11}$, i.e. $s_{11}\in R_{11}^*$.

 Since $s\cdot \begin{pmatrix}
 0& 0\\ s_{21} & 0
 \end{pmatrix} =0$, the element $\begin{pmatrix}
 0& 0\\ s_{21} & 0
 \end{pmatrix} $ belongs to the ideal $\ga$. Then $s_{11}= \begin{pmatrix}
 s_{11}& 0\\ 0 & 0
 \end{pmatrix} \in S$ since $S+\ga \subseteq S$, and $\ker (\cdot s_{11}) = \begin{pmatrix}
 0& R_{12}\\ 0 & R_{22}
 \end{pmatrix}$ since $s_{11}\in R_{11}^*$. Since $s_{11}R_{21}=0$, we have the inclusion
 $$ \begin{pmatrix}
 0& R_{12}\\ R_{21} & R_{22}
 \end{pmatrix}\subseteq \ga .$$
 Hence, $R_{12}R_{21}\subseteq \ga$.

 (i) $R_{12}R_{21}=0$: Suppose that $R_{12}R_{21}\neq 0$, we seek a contradiction. Fix a nonzero element $a=\begin{pmatrix}
 a& 0\\ 0 & 0
 \end{pmatrix}\in R_{12}R_{21}$. Then,  for some element
$ t= \begin{pmatrix}
 t_{11}& t_{12}\\ t_{21} & t_{22}
 \end{pmatrix}\in S$,
 $$ 0=ta = \begin{pmatrix}
 t_{11}& t_{12}\\ t_{21} & t_{22}
 \end{pmatrix} \begin{pmatrix}
 a& 0\\ 0 & 0
 \end{pmatrix} = \begin{pmatrix}
 t_{11}a& 0\\ t_{21}a & 0
 \end{pmatrix}.$$  In particular, $t_{11}a=0$. Since $\begin{pmatrix}
 0& t_{12}\\ t_{21} & t_{22}
 \end{pmatrix}\in \ga $, the element $t_{11}\in S$. The kernel $\ker_{R_{11}}(\cdot t_{11})$ of the map $\cdot t_{11}:R_{11}\ra R_{11}$ is  nonzero (otherwise $t_{11}\in R_{11}^*$ since $R_{11}$  is a left Artinian ring but  $t_{11}a=0$ and $0\neq a\in R_{11}$, a contradiction). Then
 $$\ker_R(\cdot t_{11})\supseteq \ker_{R_{11}}(\cdot t_{11})\bigoplus \begin{pmatrix}
 0& R_{12}\\ 0 & R_{22}
 \end{pmatrix}\supsetneqq \begin{pmatrix}
 0& R_{12}\\ 0 & R_{22}
 \end{pmatrix}=\ker_R(\cdot s).$$
This contradicts to the maximality of $\ker_R(\cdot s)$. Therefore, $R_{12}R_{21}=0$.

(ii) $S\subseteq \begin{pmatrix}
 R_{11}^*& R_{12}\\ R_{21} & R_{22}
 \end{pmatrix}$: Suppose that this inclusion does not hold, i.e. there exists an element $s'=\begin{pmatrix}
 s_{11}'& s_{12}'\\ s_{21}' & s_{22}'
 \end{pmatrix}\in S$ with $s_{11}'\not\in R_{11}^*$, we seek a contradiction. Then $s_{11}'\in S$ since $S+\ga \subseteq S$ and
 $$\begin{pmatrix}
 0& s_{12}'\\ s_{21}' & s_{22}'
 \end{pmatrix}\in\begin{pmatrix}
 0& R_{12}\\ R_{21} & R_{22}
 \end{pmatrix}\subseteq \ga .$$ By Lemma \ref{a9Feb13}.(3), $\ker_{R_{11}}(\cdot s_{11}')\neq 0$ (since $s_{11}'\not\in R_{11}^*$). Then
$$\ker_R(\cdot s_{11}')\supseteq \ker_{R_{11}}(\cdot s_{11}')\bigoplus \begin{pmatrix}
 0& R_{12}\\ 0 & R_{22}
 \end{pmatrix}\supsetneqq \begin{pmatrix}
 0& R_{12}\\ 0 & R_{22}
 \end{pmatrix}=\ker_R(\cdot s).$$
 This contradicts to the maximality of $\ker_R(\cdot s)$.

 (iii) $R_{12}=0$: Let $a_{12}= \begin{pmatrix}
0 & a_{12}\\ 0 & 0
 \end{pmatrix}\in R_{12}$. Since $R_{12}\subseteq \ga$, there is an element $t=\begin{pmatrix}
 t_{11}& *\\ * & *
 \end{pmatrix}\in S$ such that $0=ta_{12} = \begin{pmatrix}
0 & t_{11}a_{12}\\ 0 & *
 \end{pmatrix}$.  We must have $a_{12}=0$ since $t_{11}\in R_{11}^*$, by (ii).

 (iv) $\ga = \begin{pmatrix}
 0&0\\ R_{21} & R_{22}
 \end{pmatrix}=(1-e_1)R$:  The second equality is obvious.   By (ii) and (iii), $\ga \subseteq \begin{pmatrix}
 0&0\\ R_{21} & R_{22}
 \end{pmatrix}$. By (\ref{s=00}),  $s\begin{pmatrix}
 0&0\\ R_{21} & R_{22}
 \end{pmatrix}=0$. Hence, $\ga = \begin{pmatrix}
 0&0\\ R_{21} & R_{22}
 \end{pmatrix}$.

(v) $S_{e_1}=\{ 1, e_1\}\in \Den_l(R, \ga )$ {\em and }  $S_{e_1}^{-1}R\simeq R/ (1-e_1)R= R/\ga \simeq S^{-1}R$ {\em are $R$-isomorphisms}: By (iii), (iv) and  Proposition \ref{c24Dec12}, $S_{e_1}\in \Den_l(R, \ga )$ and $S_{e_1}^{-1}R\simeq R/ \ga$ (an $R$-isomorphism). By Lemma \ref{a9Feb13}.(2), $S^{-1}R\simeq R/ \ga$, an $R$-isomorphism.

2.  (a) Trivial.

(c) The first two equalities   of the statement (c) follow  from statement 1 and Proposition \ref{c24Dec12}, the last equality follows from (ii).

(b)  Suppose that $\ga\neq 0$.  Let $e=e_1$ and $e_2=1-e_1$  be as above.
By Lemma \ref{c24Dec12}, $\ga = (1-e_1)R$, hence $\ga \not\subseteq \rad (R)$ (otherwise, $0\neq e_2=1-e_1\in \rad (R)$, a contradiction).  Notice that $\rad (R)= \begin{pmatrix}
 \rad (R_{11})&0\\ R_{21}& \rad (R_{22})
 \end{pmatrix}$ and
 $$ (\ga +\rad (R))/\rad (R)= \begin{pmatrix}
 0&0\\ 0& R_{22}/\rad (R_{22})
 \end{pmatrix}\subset \begin{pmatrix}
 R_{11}/\rad (R_{11})&0\\ 0& R_{22}/\rad (R_{22})
 \end{pmatrix}=R/\rad (R).$$
Therefore, $\overline{e}_2=e_2+\rad (R) = 1_{I(\ga )}+\rad (R)$ and $\overline{e}_1=e_1+\rad (R) =1- 1_{I(\ga )}+\rad (R)=1_{CI(\ga )}+\rad (R)$. Hence, $e_2= u1_{I(\ga )}u^{-1}$ and $e_1=u1_{CI(\ga )}u^{-1}$ for some unit $u\in R^*$.  Then $\ga = u^{-1} \ga u = u^{-1} e_2Ru= 1_{ I(\ga )} R$.  $\Box $

An idempotent $e$ of a ring $R$ is called a {\em left triangular idempotent} if $eR(1-e)=0$. So, an idempotent $e$ of $R$ is a left triangular iff the ring $R=\begin{pmatrix}
R _{11}& 0\\ R_{21} & R_{22}
 \end{pmatrix}$ is {\em left triangular} iff $S_e=\{ 1,e\}\in \Den_l(R)$, by Proposition \ref{c24Dec12} iff $e$ is a left denominator idempotent of $R$.

 Let $R$ be a left Artinian ring with (\ref{1=S1is}). Let
\begin{equation}\label{CIp=R}
 \CI_l':= \CI_l'(R):=\{ e_I:=\sum_{i\in I} 1_i\; | \; e_IR(1-e_I)=0\}\;\; {\rm and}\;\; \CI''_l(R):=\CI'_l(R)\backslash \{ 1\}
\end{equation}
 where $I$ is a {\em nonempty} subset of $\{ 1, \ldots , s\}$. Notice that $1-e_I=e_{CI}$ where $CI:=\{1, \ldots , s\}\backslash I$ and $e_\emptyset :=0$.  Clearly, $1\in \CI_l'$. The set $\CI_l'$ is a finite nonempty set that consists of left triangular idempotents of the ring $R$. The set $\CI_l'(R)$ is one of the most important objects as far as left denominator sets and left localizations of $R$ are concerned (Theorem \ref{A11Feb13}).

Let, for a moment, $R$ be an arbitrary ring.
 For a ring $A\in \Loc_l(R)$, let $[A]$ be the isomorphism class of $A$. We usually drop the brackets. Let  $\Loc_l(R)/\simeq$ be the set of isomorphism classes of left localizations of the ring $R$.
  The groups $\Aut (R)$ and $\Inn (R)$ act on the sets $\Den_l(R)$ and $\Loc_l(R)$ in the obvious way: the action of an automorphism $\s \in \Aut (R)$ on $S\in \Den_l(R)$ and $S^{-1}R$ is defined as $\s (S)$ and $\s (S)^{-1} \s (R)= \s (S)^{-1} R$. Let $\Den_l(R)/G$ and $\Loc_l(R)/G$ be the sets of $G$-orbits of the groups $G= \Aut (R)$, $\Inn (R)$ in $\Den_l(R)$ and $\Loc_l(R)$, respectively. Every $\Aut (R)$-orbit is a union of $\Inn (R)$-orbits. Every isomorphism class $[A]\in \Loc_l(R)/\simeq $ is a union of $\Aut (R)$-orbits and of $\Inn (R)$-orbits.

$\noindent $

By Lemma \ref{a9Feb13}.(4), for a left Artinian ring $R$ the group $\Inn (R)$ acts {\em trivially} on $\Loc_l(R)$, i.e. each element of $\Loc_l(R)$ is an $\Inn (R)$-orbit. Recall that the set $\Loc_l(R)$ consists of $R$-isomorphism classes of left localizations of the ring $R$ and if two left localizations $S^{-1}R$ and $S'^{-1}R$ are $R$-isomorphic we write $S^{-1}R=S'^{-1}R$.

$\noindent $

 The following theorem shows that, for a left Artinian ring $R$, up to isomorphism, there are only finitely many left localizations. Moreover, there are only finitely many left localizations up to $R$-automorphism, i.e. the set $\Loc_l(R)$ is finite. The set $\Ass_l(R)$ is explicitly described and it is also a finite set.

\begin{theorem}\label{A11Feb13}
Let $R$ be a left Artinian ring. Then
\begin{enumerate}
\item The map $\CI_l'(R)\ra \{ \ass (S_e)\, | \, e\in \CI_l'(R)\}$, \; $e\mapsto \ass (S_e)=(1-e)R$ is a bijection.
\item The map $\CI_l'(R)\ra \CI_l(R)/\sim $, $e\mapsto [e]$ is a bijection.
\item The map $\CI_l'(R)\ra \Loc_l(R)$, $e\mapsto S_e^{-1}R=R/(1-e)R$, is a bijection. So, $\Loc_l(R)=\{ S_e^{-1}R\, | \, e\in \CI_l'(R)\}$, $|\Loc_l(R)|=|\CI_l'(R)|<\infty$ and  up to isomorphism   there are only finitely many  left localizations of the ring $R$.
    \item Let $e,f\in \CI_l'(R)$. Then $S_e^{-1}R= S_f^{-1}R$ iff $\ass (S_e) = \ass (S_f)$ iff $e=f$.
        \item Let $e,f\in \CI_l'(R)$. Then $e\geq f$ iff $[e]\geq [f]$ iff the map $S_e^{-1}R\ra S_f^{-1}R$, $e^ir\mapsto e^ir=\frac{r}{1}$, is well-defined iff $\ass (S_e) \subseteq  \ass (S_f)$.
            \item The map  $\CI_l'(R)\ra \Ass_l(R)$, $e\mapsto (1-e)R$, is a bijection, i.e. $\Ass_l(R)=\{ (1-e)R\, | \, e\in \CI_l'(R)\}$ is a finite set,  $|\Ass_l(R)|=|\CI_l'(R)|<\infty$.
\end{enumerate}
\end{theorem}

{\it Proof}. 1. By the very definition, if the idempotents $e_1$ and $e_2$ are distinct elements of the set $\CI_l'(R)$ then $\ass (S_{e_1}) \neq  \ass (S_{e_2})$.

2. Statement 2 follows from Theorem \ref{A9Feb13}, Lemma \ref{b11Feb13} and statement 1.

3. Statement 3 follows from Theorem \ref{A9Feb13} and statement 1.

4. Statement 4 follows from statements 1 and 3.

5. It is obvious that $e\geq f$ iff $[e]\geq [f]$. By Lemma \ref{a11Feb13}, $e\geq f$ iff the map $S_e^{-1}R\ra S_f^{-1}R$, $e^ir\mapsto e^ir$, is well-defined. By Lemma \ref{a11Feb13}, $e\geq f$  iff $\ass (S_e) \subseteq \ass (S_f)$.

6. Statement 6 follows from statement 1 and Theorem \ref{A9Feb13}.(1). $\Box $

$\noindent $

{\it Remark}.  Theorem \ref{A11Feb13} shows that not every ideal of a left Artinian ring $R$ belongs to $\Ass_l(R)$.

$\noindent $

{\bf The core of a left denominator set of a left Artinian ring}.

$\noindent $

 {\it Definition}, \cite{Bav-Crit-S-Simp-lQuot}. Let $R$ be a ring and $S\in \Ore_l(R)$. The {\em core} $S_c$ of the left Ore set $S$ is the set of all the elements $s\in S$ such that $\ker (s\cdot) = \ass (S)$ where $s\cdot : R\ra R$, $r\mapsto sr$.  If $S_c\neq \emptyset$ then $SS_c\subseteq S_c$.

 $\noindent $

 The next theorem is an explicit description of the core of a left denominator set of a left Artinian ring. In particular, it is a non-empty set.

\begin{theorem}\label{A15Feb14}
Let $R$ be a left Artinian ring, $S\in \Den_l(R, \ga )$ and $\ga \neq 0$. We keep the notation of Theorem \ref{A9Feb13}. Then
$S_c=\{ s\in S\, | \, (1-e)s(1-e)=0\}\neq \emptyset$, i.e.
$S_c= \{ s= \begin{pmatrix}
 s_{11}&0\\ s_{21}& 0
 \end{pmatrix}\in S\}$, see (\ref{R=Rij1}).
\end{theorem}

{\it Proof}. Let $s=\begin{pmatrix}
 s_{11}&0\\ s_{21}& s_{22}
 \end{pmatrix}\in S$. Then $s_{22}\in \ga$ and so $ s's_{22}=0$ for some element $ s'= \begin{pmatrix}
 s_{11}'& 0\\  s_{21}'&s_{22}'
 \end{pmatrix}\in S$. Then
 $$ S\ni s's= \begin{pmatrix}
 s_{11}'&0\\s_{21}'&s_{21}'
 \end{pmatrix}\begin{pmatrix}
 s_{11}&0 \\s_{21}& 0
 \end{pmatrix}+s's_{22}'=\begin{pmatrix}
 s_{11}'s_{11}&0\\ *& 0
 \end{pmatrix}.$$
 Therefore, $S_c':=\{ s\in S\, | \, (1-e)s(1-e)=0\}\neq \emptyset$ as $s's\in S$. Clearly, $S_c'\subseteq S_c$ since $\ga = \begin{pmatrix}
 0&0 \\ * & *
 \end{pmatrix}$. Conversely, if $s= \begin{pmatrix}
 s_{11}&0\\ s_{21}& s_{22}
 \end{pmatrix}\in S_c$ then $0=s\begin{pmatrix}
 0&0\\ 0& 1
 \end{pmatrix}= \begin{pmatrix}
 0&0\\ 0& s_{22}
 \end{pmatrix}$, i.e $s\in S_c'$. Therefore, $S_c= S_c'$. $\Box $


$\noindent $

{\bf The maximal left denominator sets of a left Artinian ring}. By Theorem \ref{A11Feb13}, the posets $(\CI'_l, \geq )$ and $(\CI_l/\sim , \geq )$ are isomorphic via the map $ e\mapsto [e]$. If $e_I, e_J\in \CI'_l$ then $e_I\geq e_J$ iff $I\supseteq J$ iff $(1-e_I)R\subseteq (1-e_J)R$ (Lemma \ref{a11Feb13}).  Let $ \min  \CI_l'$ be the set of minimal elements of the poset $\CI_l'$.

\begin{lemma}\label{d11Feb13}
Let $R$ be a left Artinian ring and $\CI_l'$ be as above.
\begin{enumerate}
\item If $e_I, e_J\in \CI_l'$ and $e_Ie_J\neq 0$ then $e_Ie_J=e_{I\cap J}=e_Je_I\in \CI_l'$.
\item If $e_I, e_J\in \CI_l'$ with $I\cap J=\emptyset$ then $e_{I\cup J}= e_I+e_J\in \CI_l'$.
\item If $e_I\in \min \CI_l'$ and $e_J\in \CI_l'$ then either $I\subseteq J$ or otherwise $I\cap J=\emptyset$.
    \item If $e_I, e_J\in \min \CI_l'$ and $e_I\neq e_J$ then $e_Ie_J=0$ and $I\cap J=\emptyset$.
        \item If $e_I, e_J\in \min \CI_l'$ and $e_I\neq e_J$ then $e_IRe_J=0$.
\end{enumerate}
\end{lemma}

{\it Proof}. 1. It is obvious that $e_Ie_J=e_{I\cap J}=e_Je_I$. In view of Proposition \ref{c24Dec12}, we have to show that the equalities $e_IR(1-e_I)=0 $ and $e_JR(1-e_J)=0 $ imply  the equality $e_Ie_JR(1-e_Ie_J)=0$:
$$e_Ie_JR(1-e_Ie_J)=e_Ie_JR(1-e_J+e_J-e_Ie_J)=e_Je_IR(1-e_I)e_J=0.$$

2. Since $I\cap J=\emptyset$, we have $e_Ie_J=e_Je_I=0$ and  $1-e_I-e_J= (1-e_I)(1-e_J)$. Then
$$ (e_I+e_J)R(1-e_I-e_J) = (e_I+e_J)R(1-e_I)(1-e_J) =0.$$
By Proposition \ref{c24Dec12}, $e_I+e_J\in \CI_l'$.

3. Suppose that $I\cap J\neq \emptyset$. Then, by statement 1, $e_{I\cap J}=e_Ie_J\in \CI_l'$, and  $e_I\geq e_{I\cap J}$. By the minimality of $e_I$, we must have the equality $e_I=e_{I\cap J}$, i.e. $I\subseteq J$.

4. Statement 4 follows from statement 3.

5. Using statement 4 and the fact that $e_IR(1-e_I)=0$, we see that $0=e_IR(1-e_I) e_J = e_IRe_J$.  $\Box $


Changing, if necessary, the order of the idempotents $1_1, \ldots , 1_s$ we may assume that $\min \CI_l'(R)=\{ e_{I_1}, \ldots , e_{I_t}\}$ where
\begin{equation}\label{Rtt1}
I_1=\{ 1, \ldots , d_1\}, \; I_2=\{ d_1+1, \ldots , d_1+d_2\},\ldots ,  \; I_t=\{ \sum_{i=1}^{t-1}d_i+1, \ldots , \sum_{i=1}^{t}d_i\}
\end{equation}
for some positive natural numbers $d_1, \ldots , d_t$. The set
$$ \{ e_1:=e_{I_1}, \ldots , e_t:= e_{I_t}, e_{t+1}=1-\sum_{i=1}^te_i\}$$
is the set of orthogonal idempotents of the ring $R$ such that $1=e_1+\cdots +e_t+e_{t+1}$ (it is possible that $e_{t+1}=0$). By Lemma \ref{d11Feb13}.(5), the ring $R$ can be seen as a matrix ring
$$ R=\bigoplus_{i,j=1}^{t+1}R_{ij}\;\; {\rm where}\;\; R_{ij}:= e_iRe_j.$$
By Lemma \ref{d11Feb13}.(5), $R_{ij}=0$ for all $i\neq j$ such that $1\leq i,j\leq t$. Hence, $R_{i,t+1}= e_iR(1-e_i-\sum_{j\neq i}e_j)=0$ for $i=1, \ldots , t$. So, the ring $R$ has the form
\begin{equation}\label{Rtt2}
R=\begin{pmatrix}
 R_{11}& 0& 0&  \ldots &0\\ 0& R_{22} & 0 &\ldots &0\\
 & & \ldots & &\\ 0&\ldots & 0 &R_{tt}&0\\
 R_{t+1, 1}& \ldots &R_{t+1, t-1}&R_{t+1, t}& R_{t+1, t+1}\\
 \end{pmatrix}
\end{equation}
\begin{proposition}\label{C12Feb13}
Let $R$ be a left Artinian ring and $\CI_l'$ be as above. Let $I_1, \ldots , I_t$ be non-empty subsets of $\{ 1, \ldots , s\}$ such that $I_i\cap I_j=\emptyset$ for all $i\neq j$.  Then $\min \CI_l'(R)=\{ e_{I_1}, \ldots , e_{I_t}\}$ iff the ring $R$ has the form (\ref{Rtt2})  and none of the rings $R_{11},\ldots , R_{tt}$ is a left triangular and, for each nonempty subset $I\subseteq \{ 1, \ldots , s\}\backslash \bigcup_{j=1}^tI_j$, $e_IR(1-e_I)\neq 0$.
\end{proposition}

{\it Proof}. $(\Rightarrow )$ Obvious.

$(\Leftarrow )$ The conditions that the ring $R$ has the form (\ref{Rtt2}) and none of the rings $R_{11}, \ldots , R_{tt}$ are left triangular show that $\min \CI_l'(R)\supseteq\{ e_{I_1}, \ldots , e_{I_s}\}$. Then the  last condition of the proposition means that $\min \CI_l'(R)=\{ e_{I_1}, \ldots , e_{I_s}\}$ (if the equality does not hold then there exists an element $e_I\in \min \CI_l'(R)\backslash \{ e_{I_1}, \ldots , e_{I_s}\}$, and so $e_IR(1-e_I)=0$. By Lemma \ref{d11Feb13}.(4), $I\subseteq \{ 1, \ldots , s\}\backslash \bigcup_{j=1}^sI_j$, but this is impossible since $e_IR(1-e_I)\neq 0$, by the assumption). $\Box $

The next theorem provides a description of the maximal left denominator sets of a left Artinian ring.
\begin{theorem}\label{B11Feb13}
Let $R$ be a left Artinian ring. Then
\begin{enumerate}
\item $\maxDen_l(R)=\{ T_e\, | \, e\in \min \CI_l'(R)\} $ where $T_e=\{ u\in R\, | \, u+(1-e)R\in (R/(1-e)R)^*\}$.
\item $|\maxDen_l(R)|\leq s$  ($s$ is  the number of isomorphism classes of left simple $R$-modules).
    \item  $|\maxDen_l(R)|= s$ iff $R$ is a semisimple ring.
\end{enumerate}
\end{theorem}

{\it Proof}. 1. By Proposition \ref{b27Nov12}, $\maxAss_l(R)=\assmaxDen_l(R)$.  By Theorem \ref{A11Feb13}.(1,5), $\assmaxDen_l(R) = \{ (1-e)R\, | \, e\in \min \CI_l'\}$ and for each element $e\in \min \CI_l'$, $S_e^{-1}R\simeq R/ (1-e)R$ is a largest left quotient ring $S^{-1}_{(1-e)R}R$ where $ S_{(1-e)R}$ is the largest left denominator set with $\ass (S_{(1-e)R})=(1-e)R$. By Theorem \ref{15Nov10}.(3), $S_{(1-e)R}= T_e$.

2. Statement 2 follows from Lemma \ref{d11Feb13}.(4).

3. Statement 3 follows from statement 1 and Lemma \ref{d11Feb13}.(5).
 $\Box $


\begin{corollary}\label{A13Feb13}
Let $R$ be a left Artinian ring, $S\in \Den_l(R, \ga )$, $e\in R$ be an idempotent such that $S_e\in \Den_l(R, \ga )$ and so  $S_e^{-1}R\simeq R/\ga \simeq S^{-1}R$ (by Theorem \ref{A9Feb13}.(1) and Proposition \ref{c24Dec12}), and $M$ be an $R$-module. Then
\begin{enumerate}
\item $\ga^2=\ga$.
\item $\ga$ is a projective right $R$-module.
\item $S^{-1}M=0$ iff $\ga M=M$  iff $(1-e)M = M$ iff $eM=0$.
\item For all $\ga, \ga'\in \Ass_l(R)$, $\ga\ga'= \ga\cap \ga'= \ga'\ga$.
\end{enumerate}
\end{corollary}

{\it Proof}. 1. We keep the notation of Theorem \ref{A9Feb13}. Then $\ga = R_{21}+R_{22}$ and $\ga^2= R_{22}R_{21}+R_{22}^2= R_{21}+R_{22}=\ga$.

2. Since $\ga = (1-e)R$, we have $R= \ga \oplus eR$, and so the ideal $\ga$ is a projective right $R$-module.

3. $S^{-1}M=0$ iff $S^{-1}_eM=0$ iff $M=\tor_{S_e}(M)$ iff $eM=0$ iff $M=(1-e)M= (1-e)RM= \ga M$.

4. By Theorem \ref{A11Feb13}.(6), $\ga = (1-e)R$ and $\ga' = (1-e')R$ for some idempotents $e,e'\in \CI_l'(R)$. The idempotents $1-e$ and $1-e'$ commute. So,
$ \ga \cap \ga' = (1-e)(1-e') R = (1-e)\ga'= (1-e)R\ga'= \ga\ga'$. $\Box$




{\bf Duality}. Let $R$ be an Artinian ring and $R^{op}$ be its opposite ring. Then $\Den_r(R)= \Den_l(R^{op})$. Using this fact we have analogous results for right denominator sets of $R$. We replace the subscript `l' by
`r' everywhere when dealing with `right' concepts. For example, $\CI'_r(R):= \CI'_l(R^{op})$ and $\CI''_r(R):= \CI''_l(R^{op})$. Recall that $\CI''_l(R)=\{ e_I\, | \, e_IRe_{CI}=0$ where $\emptyset \neq I\subsetneqq \{ 1, \ldots , s\}\}$. Similarly, $\CI''_r(R)=\{ e_J \, | \, e_{CJ}Re_{J}=0$ where $\emptyset \neq J\subsetneqq \{ 1, \ldots , s\}\}$. Then the map

\begin{equation}\label{DuIlr}
\CI''_l(R)\ra \CI''_r(R), \;\; e_I\mapsto 1-e_I= e_{CI},
\end{equation}
is {\em an order-reversing bijection of posets}.  Let $\min\CI''_*(R)$ and $\max\CI''_*(R)$ be the sets of minimal and maximal elements of the poset $\CI''_*(R)$ where $*= l,r$. The map (\ref{DuIlr}) induces the bijections
\begin{equation}\label{DuIlr1}
\min\CI''_l(R)\ra \max \CI''_r(R), \;\; e_I\mapsto 1-e_I= e_{CI},
\end{equation}
\begin{equation}\label{DuIlr2}
\max\CI''_l(R)\ra \min\CI''_r(R), \;\; e_I\mapsto 1-e_I= e_{CI}.
\end{equation}
Clearly,
\begin{eqnarray*}
 \Loc_l(R)\backslash [R]&=& \{  [ R/ (1-e_I)R]=  [ R/ e_{CI}R] \, | \, I\in \CI''_l(R)\} = \{  [ R/ e_JR] \, | \, J\in \CI''_r(R)\}, \\
  \Loc_r(R)\backslash [R]&=& \{  [ R/R (1-e_J)]=  [ R/ Re_{CJ}] \, | \, J\in \CI''_r(R)\} = \{  [ R/ Re_I] \, | \, I\in \CI''_l(R)\}.
\end{eqnarray*}
These equalities and the bijections (\ref{DuIlr1}) and (\ref{DuIlr2}) imply the next theorem.
\begin{theorem}\label{25Apr14}
Let $R$ be an Artinian ring. Then the map
$$\Loc_l(R)\backslash [R]\ra \Loc_r(R)\backslash [R], \;\; [ R/ (1-e_I)R]\ra [ R/ R(1-e_{CI})], $$ is an anti-isomorphism of posets (i.e. an order reversing bijection). In particular, $|\Loc_l(R)|= |\Loc_r(R)|$.
\end{theorem}


{\bf The left localization radical $\gll_R$ of a left Artinian ring $R$}. The next theorem gives an explicit description of the left localization radical $\gll = \gll_R$  of a left Artinian ring $R$ and  a criterion for $\gll \neq 0$.

\begin{theorem}\label{D12Feb13}
Let $R$ be a left Artinian ring. Then
\begin{enumerate}
\item $\gll = \bigcap_{e\in \min \CI'_l(R)}(1-e)R=\prod_{e\in \min \CI'_l(R)}(1-e)R=(1-\sum_{e\in \min  \CI'_l(R)} e)R$. If, in addition, $R$ is a right Artinian ring then $R=\prod_{f_j\in \max \CI_r''(R)}f_jR$.
\item $\gll = \begin{cases}
\oplus_{i=1}^{t+1}R_{t+1,i}\neq 0& \text{if }1\neq \sum_{e\in \min \CI_l'(R)}e,\\
0 & \text{otherwise}.\\
\end{cases}$
\item $\gll^2= \gll$ and $\gll$ is a projective right $R$-module.
\item $\gll\subseteq \rad (R)$ iff $\gll =0$.
\end{enumerate}
\end{theorem}

{\it Proof}. 1. The first equality follows from Theorem \ref{B11Feb13}.(1). The minimal idempotents $e\in \min \CI'_l(R)$ are orthogonal idempotents. Hence, $\cap_{e\in\min \CI'_l(R)}(1-e) R = (\prod_{e\in\min \CI'_l(R)}(1-e) ) R= (1-\sum_{e\in \min  \CI'_l(R)} e)R$. The last equality in statement 1 follows from (\ref{DuIlr1}).

2. Statement 2 follows from statement 1 and (\ref{Rtt2}) since $e_{t+1} = 1-\sum_{e\in \min  \CI'_l(R)} e$.

3. Statement 3 follows from statement 1.

4. The radical $\rad (R)$ is a nilpotent ideal. So, statement 4 follows from statement 3.  $\Box $

\begin{corollary}\label{a16Mar14}
 Let $R$ be a left Artinian ring and $e:= \sum_{e'\in \min \, \CI'_l(R)}e'$. Then
\begin{enumerate}
\item $S_e\in \Den_l(R)$.
\item $\ass (S_e) = \gll_R$.
\item $e$ is the least upper bound of the set $\min \, \CI_l'(R)$ in $\CI_l'(R)$.
\end{enumerate}
\end{corollary}

{\it Proof}. 1. By Lemma \ref{d11Feb13}.(4), $e'e''=0$ for all distinct idempotents $e'$, $e''\in \min \, \CI'_l(R)$. Hence, $1-e=\prod_{e'\in \min \, \CI'_l(R)}(1-e')$ and so $ eR(1-e)= (\sum_{e'\in \min \, \CI'_l(R)} e')\cdot R\cdot \prod_{e''\in \min \, \CI'_l(R)}(1-e'')=0$. Hence, $S_e\in \Den_l(R)$, by Proposition \ref{c24Dec12}.

2. By statement 1 and Proposition \ref{c24Dec12}, $\ass (S_e) = (1-e)R = \gll_R$ (Theorem \ref{D12Feb13}.(1)).

3. Clearly, $e\geq e'$ for all $e'\in \min \, \CI_l'(R)$. Given $f\in \CI'_l(R)$ such that $f\geq e'$  for all $e'\in \min \, \CI_l'(R)$, then $f\geq e$, by Lemma \ref{d11Feb13}.(3,4), and statement 3 follows.
 $\Box $



\section{Structure of left Artinian rings with zero left localization radical}\label{STLARZLR}

In this section, a characterization of left Artinian rings with zero left localization radical is given (Theorem \ref{A16Feb14}). A criterion is given for a left Artinian ring to be a left localization maximal ring (Theorem \ref{bb11Feb13}). For an Artinian ring, it is shown that if the left localization radical is zero then so is the right localization radical, and vice versa (Theorem \ref{c16Feb14}).

{\bf The left localization maximal rings}. These are precisely the rings  in which  we cannot invert anything on the left.

 $\noindent $

 {\it Definition}, \cite{larglquot}. A ring $A$ is
called a {\em left localization maximal ring} if $A= Q_l(A)$ and
$\Ass_l(A) = \{ 0\}$. A ring $A$ is called a {\em right
localization maximal ring} if $A= Q_r(A)$ and $\Ass_r(A) = \{
0\}$. 
  A ring $A$ is called a {\em localization maximal ring} (or a {\em two-sided localization maximal ring}) if $A= Q_{l,r}(A)$ and $\Ass_{l,r}(A) = \{ 0\}$.

$\noindent $


{\it Example}. Let $A$ be a simple ring. Then $Q_l(A)$ is a left
localization maximal  ring and $Q_r(A)$ is a right localization
maximal ring. In particular, a division ring is a (left; right; and two-sided) localization maximal ring. More generally, a simple  Artinian ring (i.e. the matrix ring over a division ring) is a (left;  right; and two-sided)
localization maximal ring.

$\noindent $

The next theorem is a criterion for   a left quotient ring of a
ring to be  a maximal left quotient ring.

\begin{theorem}\label{21Nov10}
\cite{larglquot} Let  a ring $A$ be a left localization of a ring
$R$, i.e. $A\in \Loc_l(R, \ga )$ for some $\ga \in \Ass_l( R)$.
Then $A\in \maxLoc_l( R)$ iff $Q_l( A) = A$ and  $\Ass_l(A) = \{
0\}$, i.e. $A$ is a left localization maximal ring.
\end{theorem}


\begin{corollary}\label{xy21Nov10}
The left localization maximal
rings are precisely the localizations of all the rings at their
maximal left denominators sets.
\end{corollary}

{\it Proof}. The statement follows from Theorem \ref{15Nov10}.(6) and
 Theorem \ref{21Nov10}.  $\Box $

$\noindent $

{\bf Criterion for a left Artinian ring to be a left localization maximal ring}. The next theorem is a criterion for a left Artinian ring to be a left localizable maximal ring.

\begin{theorem}\label{bb11Feb13}
Let $R$ be a left Artinian ring. Then the following statement are equivalent.
\begin{enumerate}
\item $R$ is a left localization maximal ring.
\item $\CI_l'(R)=\{ 1\}$.
\item Either $s=1$ or, otherwise, for every proper subset $I$ of $\{ 1, \ldots , s\}$, $ e_IR(1-e_I)\neq 0$.
\end{enumerate}
\end{theorem}

{\it Proof}. $(1 \Leftrightarrow 2)$  Theorem \ref{A11Feb13}.(3).

$(2 \Leftrightarrow 3)$  Proposition \ref{c24Dec12}. $\Box$


Theorem \ref{bb11Feb13} shows that `generically' every left Artinian ring is left localization maximal.

$\noindent $

{\bf Left-right symmetry of localization maximality for Artinian rings}.
 For Artinian  rings, the concept of `localization maximality' is left-right symmetric as the next theorem shows.

\begin{theorem}\label{B16Feb14}
An Artinian ring is left localization maximal iff it is right localization maximal.
\end{theorem}

{\it Proof}. The result is obvious if $s=1$. If $s>1$ then, by Theorem \ref{bb11Feb13}, $R$ is a left localization maximal iff $e_IRe_{CI}\neq 0$ for all proper subsets $I$ of $\{ 1, \ldots , s\}$ iff $e_{CI}Re_I\neq 0$ for all proper subsets  $I$ of $\{ 1, \ldots , s\}$ iff $R$ is a right localization maximal ring.  $\Box $

$\noindent $

{\bf Structure of left Artinian rings with zero left localization radical}.
\begin{theorem}\label{A16Feb14}
Let $R$ be a left Artinian ring and $\gll$ be its left localization radical. Then
\begin{enumerate}
\item The ring $R$ is a direct product of left localization maximal rings iff $\gll =0$ iff $\gll$ is a nilpotent ideal of $R$.
\item If $R$ is a direct product of left localization maximal rings then the direct product is unique up to order, i.e. if $R=\prod_{i=1}^tA_i = \prod_{j=1}^{t'}A_j'$ are two such direct products then $t=t'$ and $A_i = A_{\s (i)}'$ for all $i=1, \ldots , t$ where $\s$ is a permutation of the set $\{ 1, \ldots , t\}$.
\end{enumerate}
\end{theorem}

{\it Proof}. 1. By Theorem \ref{D12Feb13}, $\gll =0$ iff $\gll $ is a nilpotent ideal. It remains to establish the first `iff'  in statement 1.

$(\Leftarrow )$ If $\gll =0$ then $R$ is a finite direct product of localization maximal rings, by Theorem
 \ref{D12Feb13}.(2), Theorem \ref{bb11Feb13} and (\ref{Rtt2}).

$(\Rightarrow )$  If $R= \prod_{i=1}^tA_i$ is a direct product of left localization maximal (necessarily left Artinian) rings $A_i$. Then $\gll_{A_i}=0$ and $\gll = \prod_{i=1}^r \gll_{A_i}=0$, by Corollary \ref{b16Feb14}.

2. Let $1= e_1+\cdots + e_t$ and  $1= e_1'+\cdots + e'_{t'}$ be the sums of central idempotents that  correspond to the direct products $R=\prod_{i=1}^tA_i$ and $R = \prod_{j=1}^{t'}A_j'$. By Theorem \ref{bb11Feb13}, for each $i=1, \ldots , t$ there is a unique $j=\s (i)$ such that $e_i = e_ie'_j$. Then, by symmetry, $t=t'$ and $e_i = e'_{\s (i)}$ for some permutation $\s $ of the set $\{ 1, \ldots , t\}$.
 $\Box $


\begin{theorem}\label{c16Feb14}
Let $R$ be an Artinian ring, $\gll$ and $\gr$ be the left and right localization radicals of $R$, respectively. Then $\gll =0$ iff $\gr =0$.
\end{theorem}

{\it Proof}.  $\gll =0$ iff $R=\prod_{i=1}^sA_i$ is a direct product of left localization maximal rings (Theorem \ref{A16Feb14}.(1)) iff $R=\prod_{i=1}^sA_i$ is a direct product of right localization maximal rings (Theorem \ref{B16Feb14}) iff $\gr =0$ (Theorem \ref{A16Feb14}).  $\Box $

$\noindent $

{\it Remark}. In general, for an Artinian ring $R$, $\gll \neq  \gr$ (Corollary \ref{f23Dec12}).


\section{Characterization of the left localization radical of a left Artinian ring}\label{CLLRLAR}
The aim of this section is to introduce the little left localization radical of a ring and to give a characterization of it and of the left localization radical of a left Artinian ring (Theorem \ref{A20Feb14}).
A right ideal $\ga$ of a ring $R$ is called an {\em idempotent right ideal} if $\ga = fR$ for some idempotent $f\in R$.  A left  ideal $\gb$ of a ring $R$ is called an {\em idempotent left ideal} if $\gb = Re$ for some idempotent $e\in R$.
 Let $\CII_l(R)$  and $\CII_r(R)$ be the sets of left and right idempotent ideals of the ring $R$, respectively, and $\CIT_l(R)$ and $\CIT_r(R)$ be the sets of left and right idempotent ideals that in addition are {\em two-sided} ideals,  respectively.

\begin{lemma}\label{a20Feb14}
Let $R$ be a ring and $f,f'$ be idempotents of $R$. Then $fR= f'R$ iff $R(1-f) = R(1-f')$.
\end{lemma}

{\it Proof}. In view of left-right symmetry it suffices to prove that the implication
$(\Rightarrow )$ holds. The equality $fR= f'R$ yields the equality $f=f'f$ ($f=f'r$ for some $r\in R$, hence $f'f=f'^2e= f'r=f$). Then
 $(1-f')(1-f)= 1-f'-f+f'f= 1-f'$, and so  $R(1-f') \subseteq R(1-f)$. By symmetry, $R(1-f) = R(1-f')$.  $\Box $


\begin{lemma}\label{c20Feb14}
Let $R$ be a ring and $f$ be an idempotent of $R$. Then $fR$ is an ideal iff $(1-f)Rf=0$ iff $1-f\in \CI_l(R)$ iff $R(1-f)$ is an ideal.
\end{lemma}

{\it Proof}. Notice that $R=fR\oplus (1-f)R$. Then $fR$ is an ideal of $R$ iff $(1-f)RfR\subseteq fR$ iff $(1-f)Rf=0$ iff $(1-f)\in \CI_l(R)$, by (\ref{IlR=d}).
 Similarly, $R= Rf\oplus R(1-f)$. So, $R(1-f)$ is an ideal iff $R(1-f)Rf\subseteq R(1-f)$ iff $(1-f)Rf=0$.
$\Box $

\begin{corollary}\label{b20Feb14}
Let $R$ be a ring.
\begin{enumerate}
\item The map $\CII_r(R)\ra \CII_l(R)$, $ fR\mapsto R(1-f)$, is  a bijection with the inverse $Re\mapsto (1-e)R$.
\item The map $\CIT_r(R)\ra \CIT_l(R)$, $ fR\mapsto R(1-f)$, is  a bijection with the inverse $Re\mapsto (1-e)R$.
\end{enumerate}
\end{corollary}

{\it Proof}. 1. By Lemma \ref{a20Feb14}, the maps $fR\mapsto R(1-f)$ and $Re\mapsto (1-e)R$ are well-defined. Clearly, they are mutually inverse.

2. Statement 2 follows from statement 1 and Lemma \ref{c20Feb14}. $\Box $



\begin{corollary}\label{d20Feb14}
Let $R$ be a ring.
\begin{enumerate}
\item The map $\CI_l(R)/\sim \,  \ra \CIT_r(R)$, $ [e]\mapsto (1-e)R$, is  a bijection with the inverse $fR\mapsto [1-f]$.
\item  $\CIT_r(R)=\{ \ass (S_e)\,  | \, S_e\in \IDen_l(R)\}$.
\end{enumerate}
\end{corollary}

{\it Proof}. 1. Statement 1 follows from Lemma \ref{b11Feb13} and Lemma \ref{c20Feb14}.

2. Statement 2 follows from statement 1. $\Box $

$\noindent $

Let $R$ be an arbitrary ring. The intersection
$$ \gll' := \gll_R':=\bigcap_{\ga \in \Ass_l(R)\backslash \{ 0\} }\ga$$ is called the {\em little left localization radical} provided $\Ass_l(R)\backslash \{ 0\} \neq \emptyset$ and $\gll':=0$, otherwise. The set $(\CIT_r(R), \leq )$ is a poset. For an arbitrary poset $(P, \leq )$, we denote by $\Max (P)$ and $\Min (P)$ the sets of maximal and minimal elements of $P$, respectively.

$\noindent $

{\bf Characterization of the (little)  left localization radical of a left Artinian ring}.

\begin{theorem}\label{A20Feb14}
Let $R$ be a left Artinian  ring and $\gll$ be its left localization radical of $R$. Then
\begin{enumerate}
\item $\gll =\bigcap_{\ga \in \Max (\CIT_r(R))\backslash R} \ga = \prod_{\ga \in \Max (\CIT_r(R))\backslash R}\ga $.
\item Suppose that $\CIT_r(R)\neq \{ 0\}$. Then $\gll' = \prod_{\ga \in \CIT_r(R)\backslash \{ 0\}}\ga =\bigcap_{\ga \in \CIT_r(R)\backslash \{ 0\}} \ga =\bigcap_{\ga \in \Min (\CIT_r(R))\backslash \{ 0\}} \ga = \prod_{\ga \in \Min (\CIT_r(R))\backslash R}\ga $.
\end{enumerate}
\end{theorem}

{\it Proof}.  Statements 1 and 2 follow from Corollary \ref{d20Feb14}.(1) and Corollary \ref{A13Feb13}.  $\Box $

$\noindent $


\section{ Description of left denominator sets of a left Artinian ring}\label{IAMMMA}

In this section, a description of left denominator sets of a left Artinian ring $R$ is given (Theorem \ref{F12Feb13}), the sets of left localizable and non-localizable   elements are described (Proposition \ref{c11Feb13} and Proposition \ref{E12Feb13}). Theorem \ref{C11Feb13} is a criterion for the powers of a non-nilpotent element to be a left denominator set. Theorem \ref{A12Feb13} describes the set $\CC_l(R)$ of completely left localizable elements of $R$. Theorem \ref{B12Feb13} is a criterion for $\CC_l(R) = R^*$.

{\bf The set of left localizable elements of a left Artinian ring}.

\begin{lemma}\label{b11Dec12}
Let $S\in \Den_l(R, \ga )$ (respectively, $S\in \Den (R, \ga )$),
$\s : R\ra S^{-1}R$, $r\mapsto \frac{r}{1}$, and $G:=(S^{-1}R)^*$
be the group of units of the ring $S^{-1}R$. Then $S':=\s^{-1}
(G)\in \Den_l(R, \ga )$ (respectively, $S':=\s^{-1} (G)\in \Den
(R, \ga )$) and $S^{-1}R\simeq S'^{-1}R$.
\end{lemma}

{\it Proof}. It suffices to prove the statement for $S\in
\Den_l(R, \ga )$. Clearly, $S'$ is a multiplicative set.

(i) $S'\in \Ore_l(R)$: We have to show that for given elements $s'\in S'$
and $r\in R$ there are elements $s''\in S'$ and $r'\in R$ such that
$s''r=r's'$.
$$ \s (r) = \s (r) \s(s')^{-1}\s (s') = \s(s_1)^{-1} \s (r_1) \s
(s')$$ for some elements $s_1\in S$ and $r_1\in R$.  Hence, $s_2(
s_1r-r_1s')=0$ for some element $s_2\in S$. It suffices to take
$s''= s_2s_1$ and $ r'=s_2r_1$.

(ii) $\ass (S')=\ga$: If $s'a=0$ for some elements $s'\in S'$ and
$a\in R$ then $\s (s') \s (a) =0$ and so $a\in \ker (\s ) = \ga$
since $\s (s') \in G$.

(iii) $S'\in \Denl(R, \ga )$: If $as'=0$ for some elements $a\in
R$ and $s'\in S'$ then $\s (a) \s (s') =0$ and so $a\in \ker ( \s
) = \ga $ since $\s (s')\in G$.

(iv) $S'^{-1}R\simeq S^{-1}R$:  By the universal property of left localizations, the map $S'^{-1}R\ra S^{-1}R$, $s^{-1}r\mapsto s^{-1}r$, is a well-defined ring homomorphism which is obviously an epimorphism. Its  kernel is equal to zero: if  $s^{-1}r=0$ then $r\in \ga = \ass (S)= \ass (S')$, and so $s^{-1}r=0$. $\Box $


$\noindent $

{\it Definition}, \cite{Bav-Crit-S-Simp-lQuot}. An element $r$ of a ring $R$ is called a {\em left localizable element} if there  exists a left denominator set $S$  of
 $R$ such that $r\in S$ (and so the element $\frac{r}{1}\neq 0$ is invertible in the ring
 $S^{-1}R$), equivalently,  if there  exists a left denominator set $T$  of
 $R$ such that the element $\frac{r}{1}$ is invertible in the ring
 $T^{-1}R$ (Lemma \ref{b11Dec12}). The set of left localizable elements is denoted $\CL_l(R)$.

$\noindent $

Clearly, 
\begin{equation}\label{CLeU}
\CL_l(R)=\bigcup_{S\in \maxDen_l(R)} S.
\end{equation}

Similarly, a {\em right localizable element} is defined and let
$\CL_r(R)$ be the set of right localizable elements of the ring
$R$. The elements of the {\em set of left and right localizable
elements},
$$ \CL_{l,r}(R) = \CL_l(R)\cap \CL_r(R),$$
 are called {\em left and right localizable elements}. An element $r\in R$ is called a {\em localizable element} if there exists a
  (left and right) denominator set $S\in \Den (R)$ such that $r\in S$, equivalently,  if there exists a
  (left and right) denominator set $T\in \Den (R)$ such that  the element $\frac{r}{1}$ is invertible in the ring
 $T^{-1}R$ (Lemma \ref{b11Dec12}).  The set of all localizable elements of the ring $R$ is denoted by $\CL (R)$. Clearly,
 $$ \CL (R) \subseteq \CL_{l,r}(R).$$
 The sets
 $$ \CN\CL_l(R) := R\backslash \CL_l (R), \;\; \CN\CL_r(R) := R\backslash \CL_r (R), \;\;
 \CN\CL_{l,r}(R) := R\backslash \CL_{l,r} (R), \;\;\CN\CL (R) := R\backslash \CL (R),$$

are called the {\em  sets of left non-localizable; of right
non-localizable; of left and right non-localizable; of
non-localizable elements}, respectively, \cite{Bav-Crit-S-Simp-lQuot}. The elements of these
sets are called correspondingly (eg, an element $r\in \CN\CL_l(R)$
is called a {\em left non-localizable element}).

The next corollary is an explicit description of the set $\CL_l(R)$ of left localizable elements for a left Artinian ring.

\begin{proposition}\label{c11Feb13}
Let $R$ be a left Artinian ring. Then $\CL_l(R) = \{ r\in R\, | \, r+(1-e)R\in  (R/ (1-e)R)^*$ for some $e\in \min \CI_l'(R)\}$.
\end{proposition}

{\it Proof}. Let $\CR$ be the RHS of the equality of the corollary. Notice that $$\CL_l(R) = \bigcup_{S\in \maxDen_l(R)} S = \bigcup_{e\in \min \CI_l'(R)} T_e=\CR,$$ by Theorem \ref{B11Feb13}.(1). $\Box $

$\noindent $

{\bf The set of left non-localizable elements of a left Artinian ring}. The next proposition describes the set of left left non-localizable elements of a left Artinian ring and gives a criterion when it is an ideal.

\begin{proposition}\label{E12Feb13}
Let $R$ be a left Artinian ring. Then
\begin{enumerate}
\item The set $\CN\CL_l(R)$ of left non-localizable elements of $R$ is equal to the set  $\{ r\in R\, | \, r+(1-e)R\not\in (R/ (1-e)R)^*$ for all $e\in \min \CI_l'(R)\}$.
\item $R\cdot \CN\CL_l(R)\cdot R\subseteq \CN\CL_l(R)$. In particular, $\CN\CL_l(R)\cdot \CN\CL_l(R)\subseteq \CN\CL_l(R)$.
\item $\CN\CL_l(R)$ is an ideal of the ring $R$ iff the ring $R/(1-e)R$ is a division ring for all $e\in \min \CI_l'(R)$.
\end{enumerate}
\end{proposition}

{\it Proof}. 1. Statement 1 follows from Proposition \ref{c11Feb13}.

2. Statement 2 follows from statement 1 and the fact that in a left Artinian ring every one-sided invertible element is a unit.

3. $(\Leftarrow )$  Trivial.

$(\Rightarrow )$ In view of statement 2, we have to show that $\CN\CL_l(R)+\CN\CL_l(R)\subseteq \CN\CL_l(R)$. Suppose that the ring $A=R/(1-e)R$ is not a division ring for some idempotent $e\in \min \CI_L'(R)$, we seek a contradiction. Then, up to order, the ring  $\bA := A/ \rad (A) = \prod_{i=1}^m\bR_i$ is a direct product of matrix rings $\bR_i = M_{n_i}(D_i)$  where $m\leq n$.  Notice that an  element $r\in R$ is a unit iff $r+\rad (R) \in \bR$ is a unit.

Case $m=1$ and  $n_1\geq 2$. Then the elements $a= E_{11}(1)$ and
 $b= \sum_{i=2}^{n_1}E_{ii}(1)$ belong to $\CN\CL_l(R)$ but their sum does not, $a+b\not\in \CN\CL_l(R)$, a contradiction.

Case $m\geq 2$. Then the elements $a=1_1$ and $b=\sum_{i=2}^m1_i$ belong to $\CN\CL_l(R)$ but their sum, $a+b$,  does not, a contradiction. $\Box $

$\noindent $

{\bf Criterion for the powers of an element to be a left denominator set}. For a left Artinian ring $R$, the following theorem is an explicit criterion for the powers of a non-nilpotent element of $R$ to be a left denominator set.

\begin{theorem}\label{C11Feb13}
Let $R$ be a left Artinian ring, $\CI_l'(R)$ be as above, $s\in R$ be a non-nilpotent element of $R$, $e=e(s)$ be the idempotent associated with the element $s$, $S_e=\{ 1,e\}$ and $S_s=\{ s^i\, | \, i\in \N\}$. The following statements are equivalent.
\begin{enumerate}
\item $S_s\in \Den_l(R)$.
\item $S_e\in \Den_l(R)$ and $(1-e)s(1-e)$ is a nilpotent element.
\item $eR(1-e)=0$ and $(1-e)s(1-e)$ is a nilpotent element.
\end{enumerate}
If one of the equivalent conditions holds then $\ass (S_s) = (1-e)R$ and $S_s^{-1}R \simeq S_e^{-1}R \simeq R/ (1-e)R$, the core $S_{s,c}$ of the left denominator set $S_s$ is equal to $\{ s^i\, | \, i\in \N , (1-e) s^i(1-e)=0\}$.
\end{theorem}

{\it Proof}. $(2\Leftrightarrow 3)$  Statements 2 and 3 are equivalent by Lemma \ref{c24Dec12}.

$(1\Rightarrow 2)$ Suppose that $S_s\in \Den_l(R)$. By Theorem \ref{A9Feb13} there exists  $S_{e_1}\in \Den_l(R)$ such that   $S_s^{-1}R\simeq S_{e_1}^{-1}R\simeq R/ (1-e_1)R$ and $\ass (S_s)=\ass (S_{e_1}) = (1-e_1)R=:\ga $. With respect to the matrix decomposition $R=\oplus_{i,j=1}^2 R_{ij}$ associated with the idempotent $e_1$, the ring $R$ is left triangular, i.e. $R_{12}=0$ and the element $s$ has the form $s=\begin{pmatrix}
 u& 0\\v & w
 \end{pmatrix}$ where $u\in R_{11}^*$ (Theorem \ref{A9Feb13}). For all $i\geq 1$, $s^i = \begin{pmatrix}
 u^i& 0\\v_i & w^i
 \end{pmatrix}$ where $v_i\in R_{21}$. Since $e_2:= 1-e_1\in \ga = e_2R$, $0=s^i e_2 = \begin{pmatrix}
 0& 0\\0 & w^i
 \end{pmatrix}$ for some $i\geq 1$, and so $w^i =0$.

 (i) $Re_1=Rs^i$ {\em for all} $i\geq 1$ {\em such that} $ w^i=0$:  Since $s^i = \begin{pmatrix}
 u^i& 0\\v_i & 0
 \end{pmatrix}\in Re_1$, we see that $Re_1\supseteq Rs^i$. Since $s^i =\begin{pmatrix}
 u^i& 0\\ 0 & 1
 \end{pmatrix}\begin{pmatrix}
 1& 0\\v_i & 0
 \end{pmatrix}$ and $u^i\in R_{11}^*$, $Rs^i = R\begin{pmatrix}
 1& 0\\v_i & 0
 \end{pmatrix}=Re_1$ since $\begin{pmatrix}
 a& 0\\b & 0
 \end{pmatrix} \begin{pmatrix}
 1& 0\\v_i & 0
 \end{pmatrix}=\begin{pmatrix}
 a& 0\\b & 0
 \end{pmatrix}$.

 (ii) {\em If} $w^i=0$ {\em then } $Rs^i = Rs^{i+1}$: Since $s^{i+1} = ss^i =\begin{pmatrix}
 u& 0\\v & w
 \end{pmatrix}\begin{pmatrix}
 u^i& 0\\v_i & 0
 \end{pmatrix} =\begin{pmatrix}
 u^{i+1}& 0\\v_{i+1} & 0
 \end{pmatrix}$, have $Rs^{i+1} = Re_1= Rs^i$, by (i).

 (iii) $Re_1= Re$: By (i) and (ii), $Re_1= Rs^i = Re$.

 (iv) $S_e\in \Den_l(R, \ga )$: By (iii) and Lemma \ref{a20Feb14}, $\ga = (1-e_1) R = (1-e)R$. Since $(1-e) R$ is an ideal of $R$, we must have $eR(1-e) \subseteq eR\cap (1-e) R=0$, i.e. $eR(1-e) =0$. By Proposition \ref{c24Dec12}, $S_e \in \Den_l(R, \ga )$.

 Clearly, $S_s^{-1}R= R/\ga \simeq S_e^{-1}R$. The ring $R$ can be seen as the matrix ring associated with the idempotent $e$, see (\ref{MRawe}), $R=\begin{pmatrix}
 R_{11}' & 0 \\ R_{21}' & R_{22}'
 \end{pmatrix}$. Then $s= \begin{pmatrix}
 u' & 0 \\ v' & w'
 \end{pmatrix}$ and $s^i=  \begin{pmatrix}
 u'^i & 0 \\ v'_i & w'^i
 \end{pmatrix}$ for all $i\geq 1$ where $v_i'\in R_{21}'$. Since $1-e\in\ga$ we must have $0=s^i (1-e) =  \begin{pmatrix}
 0 & 0 \\ 0 & w'^i
 \end{pmatrix}$  for some $i\geq 1$, i.e. $w'= (1-e)s(1-e)$ is a nilpotent element.

 $(1\Leftarrow 2)$ This implication and all the statements of the theorem is a particular case of Lemma \ref{dd11Feb13}. $\Box $


\begin{lemma}\label{dd11Feb13}
Let $R$ be a left triangular matrix ring $R= \begin{pmatrix}
 R_{11}& 0\\ R_{21} & R_{22}
 \end{pmatrix}$ where $R_{11}$ and $R_{22}$ are arbitrary rings and $R_{21}$ is an arbitrary $(R_{22}, R_{11})$-bimodule. Let $s= \begin{pmatrix}
 u& 0\\ v & w
 \end{pmatrix}\in R$ where $u\in R_{11}^*$ and $w$ a nilpotent element of the ring $R_{22}$. Then the set $S_s=\{ s^i\, | \, i\in \N\}$ is a left denominator set of the ring $R$ such that $\ass (S_s) = \begin{pmatrix}
 0& 0\\ R_{21} & R_{22}
 \end{pmatrix}$, $S_s^{-1}R \simeq R/ \ass (S_s) \simeq R_{11}$ and the core $S_{s,c}$ of $S_s$ is equal to $\{ s^i \, | \, i\in \N, w^i=0\}$.
\end{lemma}

{\it Proof}. First, we prove all the statements but the one about the core. Then,  the statement about the core will follow, see (v). Fix a natural number $j$ such that $w^j=0$, and so $s^j = \begin{pmatrix}
 u^j& 0\\ * & 0
 \end{pmatrix}$. Notice that for each natural number $k\geq 1$,  $S_s\in \Den_l(R)$ iff $S_{s^k}\in \Den_l(R)$; and if $S_s\in \Den_l(R)$ then  $\ass (S_s)= \ass (S_{s^k})$ and $S_s^{-1}R\simeq S_{s^k}^{-1}R$. We may assume that $w=0$ from the very beginning, i.e. $s=\begin{pmatrix}
 u& 0\\ v & 0
 \end{pmatrix}$. Notice that $t=\begin{pmatrix}
 1& 0\\ -vu^{-1} & 1
 \end{pmatrix}\in R^*$ and
\begin{equation}\label{wts}
\o_t(s)= \begin{pmatrix}
 1& 0\\-vu^{-1}& 1
 \end{pmatrix}\begin{pmatrix}
 u& 0\\v& 0
 \end{pmatrix} \begin{pmatrix}
 1& 0\\-vu^{-1}& 1
 \end{pmatrix}^{-1}=\begin{pmatrix}
 u& 0\\ 0& 0
 \end{pmatrix}.
 \end{equation}
Every ideal of the ring $R$ is invariant under the inner automorphisms. So, by replacing the element $s$ by $\o_t(s)$ we may assume that $s=\begin{pmatrix}
 u& 0\\ 0 & 0
 \end{pmatrix}$.

(i) $S_s\in \Ore_l(R)$: We have to show that for given elements $s^i\in S_s$ and $r= \begin{pmatrix}
 a& 0\\ b & c
 \end{pmatrix}\in R$, there are elements $s^j\in S_s$ and $r'= \begin{pmatrix}
 a'& 0\\ b' & c'
 \end{pmatrix}\in R$ such that $s^jr=r's^i$. Notice that
 $$ sr = \begin{pmatrix}
 ua& 0\\ 0 & 0
 \end{pmatrix}= \begin{pmatrix}
 uau^{-i}& 0\\ 0 & 0
 \end{pmatrix}s^i.$$ So, it suffices to take $j=1$.

 (ii) $\ass (S_s) = \ker (s\cdot ) =\ga$ where $\ga :=  \begin{pmatrix}
 0& 0\\ R_{21} & R_{22}
 \end{pmatrix}$ since $u\in R_{11}^*$.

 (iii) $S_s\in \Den_l(R, \ga )$: $\ker (\cdot s) = \begin{pmatrix}
 0& 0\\ 0 & R_{22}
 \end{pmatrix}\subseteq \ga$.

 (iv) $S_s^{-1}R \simeq R/ \ga \simeq R_{11}$ since $u\in R_{11}^*$.

 (v) Coming back to the most general situation, then $s^i = \begin{pmatrix}
 u^i& 0\\ * & w^i
 \end{pmatrix}\in S_{s,c}$ iff $w^iR_{22}=0$ iff $w^i=0$.  $\Box $

$\noindent $

{\bf Classification of all the denominator sets of a left Artinian ring}. The next theorem gives  a criterion for a multiplicative set of a left Artinian ring to be a left denominator set, also it is an explicit  description of all the left denominator sets of $R$.

\begin{theorem}\label{F12Feb13}
Let $R$ be a left Artinian ring and $S$ be a multiplicative set of $R$. The following statements are equivalent.
\begin{enumerate}
\item $S\in \Den_l(R)$.
\item There is a nonzero idempotent $e\in R$ such that $eR(1-e)=0$, $S\subseteq \begin{pmatrix}
 R_{11}^* & 0\\ R_{21} & R_{22}
 \end{pmatrix}$ and there is an element $s\in S$ such that $s= \begin{pmatrix}
 u& 0\\ v & 0
 \end{pmatrix}$ where $R=\begin{pmatrix}
 R_{11}& 0\\ R_{21} & R_{22}
 \end{pmatrix}$ is the matrix ring associated with the idempotent $e$.
\item There is a unit $\l \in R^*$, an  idempotent $e\in \CI'_l(R)$ and an element $s\in S$ such that $\l S\l^{-1} \subseteq \begin{pmatrix}
 R_{11}^*& 0\\ R_{21} & R_{22}
 \end{pmatrix}$ and $\l s\l^{-1} = \begin{pmatrix}
 u& 0\\ v & 0
 \end{pmatrix}$ where $R=\begin{pmatrix}
 R_{11}& 0\\ R_{21} & R_{22}
 \end{pmatrix}$ is the matrix ring associated with the idempotent $e$.
 \item There is an element $s\in S$ such that $S_s\in \Den_l(R)$ and the images of all the elements of $S$ in the ring $S_s^{-1}R$ are units.
\end{enumerate}
If one of the equivalent conditions holds then $\ass (S) = \ker_R(s\cdot  )$ in   cases 2 and 3 regardless of the choice of $s$.
\end{theorem}

{\it Proof}. $(1\Rightarrow 2)$  The existence of the idempotent $e$ follows from
Theorem \ref{A9Feb13}.(2). Since $e_2:= 1-e\in a\ass (S)$, there exists an element $s= \begin{pmatrix}
 u& 0\\ v & w
 \end{pmatrix}\in S$  such that $0=se_2= \begin{pmatrix}
 0& 0\\ 0 & w
 \end{pmatrix}$, i.e. $w=0$.

$(2\Leftrightarrow 3)$ This implication is obvious due to the fact that every idempotent $e\in R$ such that $eR(1-e)=0$ is conjugate to an idempotent of the set $\CI_e(R)'$.

$(2\Rightarrow 4)$ Using the inner automorphism $\o_t$ as in (\ref{wts}), and replacing the element $s$ by $\o_t(s)$ we may assume that $s= \begin{pmatrix}
 u& 0\\ 0 & 0
 \end{pmatrix}$, i.e. $v=0$ in statement 2. By Theorem \ref{C11Feb13}, $S_s\in \Den_l(R)$ and $S_s^{-1}R\simeq     R_{11}$. So, the images of all the elements of $S$ in the ring $S_s^{-1}R$ are units.

$(4\Rightarrow 1)$ This implication is a particular case of Proposition \ref{G12Feb13}.(2) where $\ga = \ass (S_s) = \cup_{i\geq 1}\ker ( s^i\cdot )$. In particular, $S\in \Den_l(R, \ga )$.  $\Box $

$\noindent $

Theorem \ref{F12Feb13} gives the following algorithm  of obtaining all the left denominator sets $S$ of a left Artinian ring $R$. Choose an idempotent $e$ of $\CI'_l(R)$, there are only finitely many of them. Then we have the triangular decomposition of the ring $R= \begin{pmatrix}
 R_{11}& 0\\ R_{21} & R_{22}
 \end{pmatrix}$ associated with the idempotent $e$. Choose an arbitrary element $s=\begin{pmatrix}
 u& 0\\ v & 0
 \end{pmatrix}$ with $u\in R_{11}^*$. Choose an arbitrary set of elements  $s_i= \begin{pmatrix}
 u_i& 0\\ v_i & w_i
 \end{pmatrix}$, $i\in I$, where $u_i\in R_{11}^*$ for all $i\in I$. Consider the monoid $S'$  generated by the set $\{ s,s_i\, | \, i\in I\}$.    Then an arbitrary  left denominator set $S$ of $R$ is of the type $\l S'\l^{-1}$ where  $\l \in R^*$.

{\bf The set of completely left localizable elements of a left Artinian ring}.

$\noindent $

{\it Definition}, \cite{Bav-Crit-S-Simp-lQuot}. For an arbitrary ring $R$, the intersection $$\CC_l(R):=\bigcap_{S\in \maxDen_l(R)}S$$
is called the set of {\em completely left localizable elements} of $R$ and an element of the set $\CC_l(R)$ is called a {\em completely localizable element}.

$\noindent $

By Proposition \ref{A8Dec12}, $S_l(R)\subseteq \CC_l(R)$. In general, this inclusion is strict, see Theorem \ref{B12Feb13}.  Moreover, Theorem \ref{B12Feb13} is a criterion for  $S_l(R)=\CC_l(R)$ for a left Artinian ring $R$ (notice that $S_l(R) = \CC_R = R^*$). The next theorem describes the set of completely left localizable elements of a left Artinian ring.

\begin{theorem}\label{A12Feb13}
Let $R$ be a left Artinian ring. Then
$\CC_l(R)=\{ r\in R\, | \, r+(1-e)R\in (R/ (1-e)R)^*$ for all $e\in \min \CI_l'(R)\}$.
\end{theorem}

{\it Proof}. The theorem follows from Theorem \ref{B11Feb13}.
$\Box $

$\noindent $

A left Artinian ring $R$ is a {\em strongly left triangular} if there is an idempotent $e_I=\sum_{i\in I}1_i$ for some proper  subset $I$ of $\{ 1, \ldots , s\}$ such that $e_IR(1-e_I)=0$. This definition does not depend on the choice of the idempotents $1_1, \ldots , 1_s$ since any two sets of them are conjugate. A left Artinian ring $R$ which is not a strongly left triangular is called a {\em non-strongly-left-triangular}.

\begin{lemma}\label{a12Feb13}
Let $R$ be a left Artinian ring.  Then
\begin{enumerate}
\item The ring $R$ is a non-strongly-left-triangular ring iff $\CI_l(R)=\{ 1\}$ iff $\CI_l'(R)=\{ 1\}$ iff $R$ is a left localization maximal ring.
\item If the ring $R$ is non-strongly-left-triangular then $\CC_l(R)=R^*$.
\end{enumerate}
\end{lemma}

{\it Proof}. 1. Theorem \ref{bb11Feb13}.

2. This follows from statement 1,  Theorem \ref{A12Feb13} and Theorem \ref{bb11Feb13} as $\CI_l(R)=\{ 1\}$.  $\Box $

$\noindent $

The next theorem is a criterion for $\CC_l(R)=R^*$.

\begin{theorem}\label{B12Feb13}
Let $R$ be a left Artinian ring.   We assume that (\ref{Rtt1}) and (\ref{Rtt2}) hold. In particular, $\min \CI_l'(R)=\{ e_{I_1}, \ldots , e_{I_t}\}$.
Then the following statements are equivalent.
\begin{enumerate}
\item $\CC_l(R)=R^*$.
\item $\bigcup_{i=1}^tI_i=\{ 1, \ldots , s\}$.
\item $R\simeq \prod_{i=1}^tR_{ii}$.
\item $R$ is a direct product $\prod_{i=1}^{t'}R_i$ of (necessarily left Artinian) non-strongly-left-triangular rings $R_i$.
    \item $\gll_R=0$.
\end{enumerate}
If one of the equivalent conditions holds then $t=t'$ and, up to order, $R_i\simeq R_{ii}$ for $i=1, \ldots , t$.
\end{theorem}

{\it Proof}. $(2\Leftrightarrow 3)$ Obvious, see (\ref{Rtt2}).

$(3\Rightarrow 4)$ Lemma \ref{a12Feb13}.

$(1\Rightarrow 2)$  Suppose that statement 2 does not hold, we seek a contradiction.  In view of the decomposition (\ref{Rtt2}), an element $r=(r_{ij})\in R=\oplus_{i,j=1}^{t+1}R_{ij}$ is a unit iff $r_{ii}\in R_{ii}^*$ for $i=1, \ldots , t+1$. By Theorem \ref{A12Feb13}, the element $a=(a_{ij})$ which is a diagonal matrix with $a_{11}=1, \ldots , a_{tt}=1$ and $ a_{t+1, t+1}=0$ is a non-zero element which is not a unit but belongs to $\CC_l(R)$, a contradiction.

$(4\Rightarrow 1)$ By Theorem \ref{c26Dec12},
$$\maxDen_l(\prod_{i=1}^{t'}R_i)=\coprod_{i=1}^{t'}\maxDen_l(R_i).$$ By Lemma \ref{a12Feb13}.(2), $\maxDen_l(R_i)=\{ R_i^*\}$. Hence, $\CC_l(R)=R^*$, by Theorem \ref{c26Dec12}. Therefore, statements 1--4 are equivalent.

 The fact that $t=t'$ follows from Lemma \ref{a12Feb13}.(1) and Lemma \ref{b12Feb13}. Since none of the rings $R_{ii}$ and $R_i$ where $i=1, \ldots , t$, is a product of rings then it is well-known and easy to prove that, up to order, $R_i\simeq R_{ii}$ (even the equality hold) for all $i$.

 $(3\Leftrightarrow 5)$ This equivalence  follows from Theorem \ref{A16Feb14}, Lemma \ref{a12Feb13}.(1) and  the  equivalence $1\Leftrightarrow 3$.
  $\Box $

\begin{lemma}\label{b12Feb13}
Let $R=\prod_{i=1}^nR_i$ be a product of left Artinian rings $R_i$. Then
\begin{enumerate}
\item $\CI_l(R)=\{ e_{i_1}+\cdots +e_{i_\nu}\, | \, e_i\in \CI_l(R_i), i, \nu =1, \ldots , n\}$.
\item $\min \CI_l(R)=\coprod_{i=1}^n \min \CI_l(R)$.
\end{enumerate}
\end{lemma}

{\it Proof}.  Straightforward.  $\Box $






$\noindent $

{\bf The rings of lower and upper triangular matrices over a division ring}. Let $D$ be a division ring, $L_n = L_n(D)$ and $U_n = U_n(D)$ be the rings of lower and upper triangular matrices, respectively, when $n\geq 2$. Let $E_{ij}$ be the matrix units where $i,j=1, \ldots , n$. The decomposition (\ref{1=S1is}) takes the form $1=1_1+\cdots +1_n$ where $1_i:= E_{ii}$.

\begin{lemma}\label{d23Dec12}
Let $n\geq 2$ and $1_{[s]}=E_{11}+\cdots +E_{ss}$. Then
\begin{enumerate}
\item $\CI_l(L_n) = \{ 1_{[s]}\, | \, [s]= \{ 1, \ldots  , s\}, s=1, \ldots , n\}$ and $\min \CI_l'(L_n) = \{ E_{11}\}$.
\item $ \maxDen_l(L_n) = \{ T_{E_{11}}\}$ and $T_{E_{11}}^{-1}L_n\simeq S_{E_{11}}^{-1}L_n\simeq D$ (see Theorem \ref{B11Feb13}.(1)).
    \item $\gll_{L_n} = \ass (T_{E_{11}})= (1-E_{11})L_n$.
\end{enumerate}
\end{lemma}

{\it Proof}. 1. Statement 1 is obvious.

2. Statement 2 follows from statement 1,  Theorem \ref{B11Feb13} and Proposition \ref{c24Dec12}.

3. Statement 3 follows from statements 1 and 2.  $\Box $


\begin{corollary}\label{e23Dec12}
Let $n\geq 2$. Then
\begin{enumerate}
\item $\CI_l(U_n) = \{ 1_{[s]'}\, | \,  [s]':= \{ s,s+1, \ldots , n\}, s=1, \ldots , n\}$ and $\min \CI_l'(U_n) = \{ E_{nn}\}$.
\item $ \maxDen_l(U_n) = \{ T_{E_{nn}}\}$ and $T^{-1}_{E_{11}}U_n\simeq S_{E_{nn}}^{-1}U_n\simeq D$.
    \item $\gll_{U_n}= (1-E_{nn})U_n$.
\end{enumerate}
\end{corollary}

{\it Proof}.  The $D$-homomorphism
$$U_n\ra L_n , \;\; E_{ij}\mapsto E_{n+1-i, n+1-j}, \;\; i=1, \ldots , n,$$ is a ring isomorphism, and the results follows  from Lemma \ref{d23Dec12}. $\Box $

$\noindent $

Let $D^{op}:= \{ d^0\, | \, d\in D\}$ be the division ring {\em opposite} to the division ring $D$, i.e. $d_1^0d_2^0:= (d_2d_1)^0$. The map $D\ra D^0$, $d\mapsto d^0$, is a ring anti-isomorphism. The map
\begin{equation}\label{TOdef}
T_0: L_n(D)\ra U_n(D^0), \;\; \sum d_{ij}E_{ij}\mapsto \sum d_{ij}^0E_{ji},
\end{equation}
is an anti-isomorphism. Let $\CI_r'(R)$ be the right analogue of $\CI_l'(R)$:
\begin{equation}\label{Crp=R}
\CI_r'(R):=\{ e_I\, | \, (1-e_I)Re_I=0\}.
\end{equation}

Let $\gr_R$ be the right localization radical of a ring $R$.
\begin{corollary}\label{f23Dec12}
Let $n\geq 2$. Then
\begin{enumerate}
\item $\CI_r(L_n) = \{ 1_{[s]'}\, | \,  [s]':= \{ s,s+1, \ldots , n\}, s=1, \ldots , n\}$ and $\min \CI_r'(L_n) = \{ E_{nn}\}$.
\item $ \maxDen_r(L_n) = \{ T'_{E_{nn}}\}$ and $L_nT'^{-1}_{E_{nn}}\simeq L_nS_{E_{nn}}^{-1}\simeq D$ where $T'_{ E_{nn}}$ is the right analogue of $T_e$, see Theorem \ref{B11Feb13}.(1).
    \item $\gr_{L_n}= L_n(1-E_{nn})$ and $ \gll_{L_n} \neq  \gr_{L_n}$.
\end{enumerate}
\end{corollary}
{\it Proof}. These follows from Corollary \ref{e23Dec12} and (\ref{TOdef}). $\Box$

\begin{corollary}\label{g23Dec12}
Let $n\geq 2$. Then
\begin{enumerate}
\item $\CI_r(U_n) = \{ 1_{[s]}\, | \, [s]= \{ 1, \ldots  , s\}, s=1, \ldots , n\}$ and $\min \CI_r'(U_n) = \{ E_{11}\}$.
\item $ \maxDen_r(R_n) = \{ T'_{E_{11}}\}$ and $U_nT'^{-1}_{E_{11}}\simeq U_nS_{E_{11}}^{-1}\simeq D$.
    \item $\gr_{U_n}= U_n(1-E_{11})$ and $\gll_{U_n}\neq  \gr_{U_n}$.
\end{enumerate}
\end{corollary}

{\it Proof}. These follows from Lemma \ref{d23Dec12} and (\ref{TOdef}). $\Box$


\section{Localizations of Artinian rings}\label{LARR}

 This section is about (left and right) denominator sets, $\Den (R)$, and localizations $\Loc (R) := \{ S^{-1}R=RS^{-1}\, | \, S\in \Den (R)\}$ of an Artinian ring $R$. The results of this section are analogous to their left versions but much more simpler due to Corollary \ref{d24Dec12}. Their proofs follow from the left analogues  in a straightforward manner and are left for the reader as an easy exercise. For the two-sided (i.e. left and right) concepts we use the same notations but the subscript `l' is dropped, eg $\maxDen (R)$ is the set of maximal (left and right) denominator sets of $R$  and $\Ass (R) =\{ \ass (S) \, | \, S\in \maxDen (R)\}$. Briefly, for (left and right) denominator sets and localizations of Artinian rings, the central idempotents play a crucial role.

Let $R$ be a left Artinian ring. It can be uniquely (up to permutation) presented as a direct product of rings
\begin{equation}\label{RRpi}
R=\prod_{i=1}^t R_i
\end{equation}
where $R_i$ are necessarily left Artinian rings none of which is a direct product of two rings. Let 
\begin{equation}\label{RRpi1}
1=\sum_{i=1}^t e_i
\end{equation}
be the corresponding sum of central idempotents and let
\begin{equation}\label{RRpi2}
\CI (R)=\{ e_i\, | \, i=1, \ldots , t\}\;\; {\rm and}\;\; \CI' (R)=\{ e_I:= \sum_{i\in I}e_i \, | \, \emptyset \neq I\subseteq \{1, \ldots , t\} \}.
\end{equation}
Clearly, $R_i=e_iR$ and none of $e_i$ is a sum of two nonzero central idempotents of $R$.

The first statement of the following theorem shows that every  localization of an Artinian ring is a (central) idempotent  localization.
\begin{theorem}\label{TA9Feb13}
Let $R$ be an  Artinian ring and $S\in \Den (R, \ga )$. Then
\begin{enumerate}
\item There exists a central  idempotent $e\in R$ such that $S_e:=\{ 1,e\}\in \Den (R, \ga )$ and the rings $S^{-1}R$ and $ S_e^{-1}R$ are $R$-isomorphic.
\item
\begin{enumerate}
\item If $\ga =0$ then $e=1$.
\item If $\ga \neq 0$ then $\ga = (1-e)R\not\subseteq \rad (R)$.
    \item 
\begin{equation}\label{TR=Rij1}
R=\begin{pmatrix}
 R_{11}& 0\\ 0 & R_{22}
 \end{pmatrix}, \;\;\; \ga =\begin{pmatrix}
 0&0 \\ 0 & R_{22}
 \end{pmatrix} \;\; {\rm and}\;\; S\subseteq \begin{pmatrix}
 R_{11}^*& 0\\ 0 & R_{22}
 \end{pmatrix},
\end{equation}
where $R_{ij}=e_iRe_j$, $e_1=e$ and $e_2=1-e_1$.
\end{enumerate}
\end{enumerate}
\end{theorem}

The following theorem shows that, for an Artinian ring $R$ there are only finitely many left localizations. Moreover, there are only finitely many left localizations up to $R$-automorphism, i.e. the set $\Loc (R)$ is finite. The set $\Ass (R)$ is explicitly described and it is also a finite set.

\begin{theorem}\label{TA11Feb13}
Let $R$ be an  Artinian ring. Then
\begin{enumerate}
\item The map $\CI' (R)\ra \{ \ass (S_e)\, | \, e\in \CI (R)\}$, \; $e\mapsto \ass (S_e)=(1-e)R$ is a bijection.
   \item The map $\CI' (R)\ra \Loc (R)$, $e\mapsto S_e^{-1}R=R/(1-e)R$, is a bijection. So, $\Loc (R)=\{ S_e^{-1}R\, | \, e\in \CI' (R)\}$, $|\Loc (R)|=|\CI' (R)|=2^t-1<\infty$ and  up to isomorphism   there are only finitely many  left localizations of the ring $R$.
   \item The map  $\CI' (R)\ra \Ass (R)$, $e\mapsto (1-e)R$, is a bijection, i.e. $\Ass (R)=\{ (1-e)R\, | \, e\in \CI' (R)\}$ is a finite set,  $|\Ass  (R)|=|\CI' (R)|=2^t-1<\infty$.
\end{enumerate}
\end{theorem}

Let $R$ be a ring and $S\in \Ore (R)$. The {\em core} $S_c$ of the  Ore set $S$ is the set of all the elements $s\in S$ such that $\ker (s\cdot)= \ker (\cdot s) = \ass (S)$ where $\cdot s : R\ra R$, $r\mapsto rs$.  If $S_c\neq \emptyset$ then $SS_cS\subseteq S_c$.
 The next theorem is an explicit description of the core of a  denominator set of an  Artinian ring. In particular, it is a non-empty set.

\begin{theorem}\label{TA15Feb14}
Let $R$ be an  Artinian ring, $S\in \Den (R, \ga )$ and $\ga \neq 0$. We keep the notation of Theorem \ref{TA9Feb13}. Then
$S_c=\{ s\in S\, | \, (1-e)s=0\}\neq \emptyset$, i.e.
$S_c= \{ s= \begin{pmatrix}
 s_{11}&0\\ 0 & 0
 \end{pmatrix}\in S\}$, see (\ref{TR=Rij1}).
\end{theorem}

\begin{theorem}\label{6Mar13}
Let $R$ be an  Artinian ring and $\CI' (R)$ be as above.  Then
\begin{enumerate}
\item  there are precisely $2^t-1$ idempotent (left and right) denominator sets. Moreover, the map
    $$ \CI' (R)  \ra \IDen (R), \;\; e_I\mapsto S_{e_I}=\{ 1, e_I\},$$
    is a bijection where $\ass (S_{e_I})= R(1-e_I)$ and $S_{e_I}^{-1}R\simeq R/ R(1-e_I)\simeq \prod_{i\in I} R_i$.
\item $\maxDen (R) =\{ S_i:= R_1\times \cdots \times R_{i-1}\times R_i^*\times R_{i+1}\times \cdots \times R_t\, | \, i=1, \ldots , t\}$ and $\maxLoc (R)= \{ R_i \, | \, i=1, \ldots , t\}$.
    \item The (left and right) localization radical $\cap_{S\in \maxDen (R)} \ass (S)$ of $R$ is equal to zero.
   \item The set $\CC (R) :=\cap_{i=1}^t S_i$ of completely (left and right) localizable elements of $R$ is the group $R^*$ of units of $R$.
       \item The set of (left and right) localizable  elements $\CL (R) =\bigcup_{i=1}^t S_i$ is equal to $\{ r=(r_1, \ldots , r_t)\in \prod_{i=1}^t R_i \, | \, r_i\in R_i^*$ for some $i\}$.
           \item The set of (left and right) non-localizable elements $\CN\CL (R) := R\backslash \CL (R)$ is equal to $\{ (r_1, \ldots , r_n) \in \prod_{i=1}^t R_i\, | \, r_i$ is a zero divisor of $R_i$ for $i=1, \ldots , t\}$.
               \item
               \begin{enumerate}
\item Every (left and right) denominator set $S\in \Den (R)$ contains precisely one central idempotent $e_I$ where  $\emptyset \neq I\in 2^{\CI (R)}$ such that $\ass (S) =R(1-e_I)$.
\item  Every (left and right) denominator set $S\in \Den (R)$ contains precisely one  central idempotent $e_I$ where  $\emptyset \neq I\in 2^{\CI (R)}$ such that $S^{-1}R\simeq S_e^{-1}R$ is an $R$-isomorphism.
\end{enumerate}
\item  Every (left and right) denominator set $S$ of $R$ is obtained in the following way: fix $e_I\in \CI'(R)$  and take a multiplicative submonoid $S$ of $(R, \cdot )$ such that for all elements $r= (r_1, \ldots , r_n) \in S$, $r_i\in R_i^*$ for $i\in I$, and there exists an element $r'=(r_1', \ldots , r_m') \in S$ with $r_j'=0$ for all $j\not\in I$. Then $S\in \Den (R)$.
\end{enumerate}
\end{theorem}

\begin{theorem}\label{A18Mar14}
Let $R$ be an Artinian ring and $r=(r_1, \ldots , r_t)\in \prod_{i=1}^tR_i$ be a non-nilpotent element. Then $S_r= \{ r_i\, | \, i\in \N \} \in \Den (R)$ iff each $r_i$ is either a unit or a nilpotent element.
\end{theorem}

\begin{corollary}\label{a9Mar13}
 All the denominator sets of $R$ consists of units iff 1 is the only central idempotent of $R$.
\end{corollary}


\section{Rings with  left Artinian left quotient ring}\label{RWLAQR}

 The aim of this section is to show that if the left quotient ring $Q_l(R)$ of a ring $R$  is a left Artinian ring then $|\maxDen_l(R)|<\infty$ (Theorem \ref{X17Feb14}.(1)). Recall that the largest left quotient ring $Q_l(R)$  of $R$ is a left Artinian ring   iff the (classical)  left quotient ring $Q_{l, cl} (R):= \CC_R^{-1}R$ is a left Artinian ring, and in this case $S_l(R) = \CC_R$, \cite{Bav-genGoldie}.

Let $R$ be a ring, $S\in \Den_l(R)$, $M$ be an $R$-module and $N,L$ be submodules of $M$. We say that $N$ and $L$ are $S$-{\em equal} in $M$  and write $N\stackrel{S}{=}L$ if $S^{-1}N= S^{-1}L$ in $S^{-1}M$. Clearly, $N$ and $L$ are $S$-equal iff, for each pair of elements $n\in N$ and $l\in L$, $sn \in L$ and $tl\in N$  for some elements $s,t\in S$.

\begin{theorem}\label{X17Feb14}
 Let $R$ be a ring such that  $Q_l(R)$ is a left Artinian ring and $s$ be the number of iso-classes of simple  left $Q_l(R)$-modules. Then
\begin{enumerate}
\item the map $\maxDen_l(R)\ra \maxDen_l(Q_l(R))$, $S\mapsto SQ_l(R)^*$, is a bijection with the inverse $T\mapsto T\cap R$. In particular, $|\maxDen_l(R)|=|\maxDen_l(Q_l(R))|\leq s<\infty$.
    \item $|\maxDen_l(R)|=s$ iff $Q_l(R)$ is a semisimple ring iff $R$ is a semiprime left Goldie ring.
\item The map $\maxAss_l(R)\ra \maxAss_l(Q_l(R))$, $\ga\mapsto S_l(R)^{-1}\ga$, is a bijection with the inverse $\gb\mapsto \gb\cap R$. In particular, $|\maxAss_l(R)|=|\maxAss_l(Q_l(R))|\leq s<\infty$.
\item For all ideals $\ga, \ga'\in \maxAss_l(R)$, $\ga\ga' \stackrel{S}{=} \ga\cap \ga'\stackrel{S}{=}\ga'\ga$ and $\ga^2 \stackrel{S}{=}\ga$.
\end{enumerate}
\end{theorem}

{\it Proof}. 1. Statement 1 follows from Proposition \ref{A8Dec12} 
  and Theorem \ref{B11Feb13}.

2. Statement 2 follows from Theorem \ref{B11Feb13}.(3).

3. Statement 2 follows from the equality $\maxAss_l(R) =\{\ass_l(S)\, | \, S\in \maxDen_l(R)\}$ (Proposition \ref{b27Nov12}), statement 1 and the fact that the largest left quotient ring $Q_l(R)$ is a left Artinian ring.

4. Statement 3 follows from  the fact  that for all ideals $\gb , \gb'\in \Ass_l(Q_l(R))$, $\gb \gb' = \gb \cap \gb' = \gb'\gb$ (Corollary \ref{A13Feb13}). $\Box $

$\noindent $

\begin{lemma}\label{x17Feb14}
Let $R$ be a ring such that $Q_l(R)$ is a left Artinian ring, $S_l= S_l(R)$ and $\Ass_l(R, S_l):=\{ \ga \in \Ass_l(R)\, | \, \ga \cap S_l=\emptyset\}$. Then for all ideals $\ga, \ga'\in \Ass_l(R)$,
\begin{enumerate}
\item $\ga^2 \stackrel{S_l}{=}\ga$.
\item $\ga\ga' \stackrel{S_l}{=} \ga\cap \ga'\stackrel{S_l}{=}\ga'\ga$
\end{enumerate}
\end{lemma}

{\it Proof}. The ring $Q_l := Q_l(R)$ is a left Artinian ring. So, for all ideals $\ga$ of $R$, $S^{-1}\ga$ are ideals in $Q_l$. Clearly, $S^{-1}\ga = Q_l$ iff $\ga \cap S_l\neq \emptyset$. So, if one of the ideals in statements 1 and 2 meets $S_l$, then statements 1 and 2 hold. We can assume that $\ga , \ga'\in \Ass_l(R, S_l)$.
 Choose $S\in \Den_l(R, \ga )$ and $S'\in \Den_l(R, \ga' )$. Let $Q_l^*$ be the group of units of the ring $Q_l$.  Then $SQ_l^* \in \Den_l( Q_l, \gb :=S_l^{-1} \ga )$ and $S'Q_l^* \in \Den_l( Q_l, \gb' :=S_l^{-1} \ga' )$. By Corollary \ref{A13Feb13},  $\gb\gb' = \gb\cap \gb' = \gb'\gb$.  Then statements 1 and 2 follow.  $\Box $

$\noindent $


$\noindent $

$${\bf Acknowledgements}$$

 The work is partly supported by  the Royal Society  and EPSRC.

\small{

Department of Pure Mathematics

University of Sheffield

Hicks Building

Sheffield S3 7RH

UK

email: v.bavula@sheffield.ac.uk}

\end{document}